\documentclass[12pt]{amsart}
\usepackage{amssymb}
\usepackage{latexsym}
\usepackage{amsmath}

\newcommand{\eenmatrix}[4]{\setlength{\unitlength}{1ex}
\raisebox{-1.5ex}[0ex][3.5ex]{\parbox{7.5ex}{\begin{picture}(5,5)(0,0)
\put(1.8,1.9){\line(0,1){4}}
\put(1.8,5.9){\line(1,0){4}}
\put(5.8,5.9){\line(0,-1){4}}
\put(5.8,1.9){\line(-1,0){4}}
\put(2.3,6.2){\parbox[b]{3ex}{\centering $#1$}}
\put(6.2,3.4){\parbox[b]{8ex}{$#2$ }}
\put(0,3.4){\parbox[b]{1.5ex}{  \raggedleft $#3$  }}
\put(2.3,0){\parbox[t]{3ex}{\centering $#4$}}
\end{picture}
}}}

\begin{document}

\title{ Actions and coactions of finite quantum groupoids on von Neumann algebras, extensions of the match pair procedure}
\author{Jean-Michel Vallin}
\address{Administrative address:  UMR CNRS 6628 Universit\'e D'Orl\'eans,
\vskip 0,3cm
  Postal address:  Institut de
Math\'ematiques de Jussieu,  Projet alg\`ebres
\vskip 0,01cm
d'op\'erateurs et repr\'esentations, Pi\`ece 7E12, 175
rue du Chevaleret 75013 Paris
\vskip 0,3cm
  e-mail: jmva @ math.jussieu.fr}


\markboth{Jean-Michel Vallin} {Actions and coactions of finite quantum groupoids on von Neumann algebras, extensions of the match pair procedure}

\begin{abstract} In this work we investigate the notion of action or coaction of a finite quantum groupoid in  von Neumann algebras context. In particular we prove a double crossed product theorem and prove the existence of an universal von Neumann algebra on which any finite groupoids acts outerly. In  previous works,  N. Andruskiewitsch  and S.Natale define for any match pair of groupoids  two $C^*$-quantum groupoids in duality, we give here an interpretation of them in terms of crossed products of groupoids using a multiplicative partial isometry which gives a  complete description of  these structures.  In a next  work  we shall give a third description of  these structures  dealing with inclusions of depth two inclusions of von Neumann algebras associated with outer actions of match pairs of groupoids ,  and a study, in the same spirit,  of an other extension of the match pair procedure.  \end{abstract}

\maketitle

\newpage
\section{Introduction}
\def\I{ I_{\mathcal H,\mathcal K}}
\def\k{\gamma^K}
\def\h{\delta^H}
\newenvironment{dm}{\hspace*{0,15in} {\it Proof:}}{\qed}

Multiplicative partial isometries ({\it mpi}) generalize Baaj and Skandalis  multiplicative unitaries in finite dimension [BS], [BBS]. They are the 
finite-dimensional version of so-called pseudo-multiplicative unitaries, which appeared first in a commutative context dealing with locally compact groupoids [Val0], and then in the general case for a very large class of 
depth two inclusions of von Neumann algebras in a common work with M.Enock [EV], who has developped the theory in the infinite dimensional  framework leading to F.Lesieur's measured quantum groupoids [L].

When it is regular, any {\it mpi} generates two involutive subalgebras of the algebra 
of all bounded linear operators on the corresponding Hilbert space. Using 
a canonical pairing, these two algebras have structures generalizing involutive Hopf algebras. The first examples of these new structures were discovered by the theoretical physicists B\"ohm, Szlachanyi and Nill [BoSz] [BoSzNi], they called them weak Hopf $C^*$-algebras.

D.Nikshych  and  L.Vainerman, using general inclusions of depth two subfactors of 
type $II_1$ with finite index  $M_0 \subset M_1$, and a 
special pairing between relative commutatnts $M'_0 \cap M_2$ and $M'_1\cap M_3$, 
gave explicit formulas for weak Hopf $C^*$-algebra structures in 
duality for these two last involutive algebras [NV2] and  found a Galois correspondence between intermediate subfactors and involutive coideals for $M'_1 \cap M_3$ [NV4].

In an algebraic construction, N.Andruskiewitsch and S.Natale in [AA], give a construction of weak Hopf $C^*$-algebras dealing with match pairs of groupoids in a sense generalizing directly the group case.

This work is an operator  algebra point of view on these match pairs of groupoids,  we use a special multiplicative partial isometry  to give an interpretation of these examples in terms of groupoids crossed products generalizing in finite dimension previous works in the quantum groups context.

In the second paragraph we give  the fondamental definitions and properties of multiplicative partial isometries and their close connection with quantum groupoids.

The third chapter is an approach of the actions of quantum groupoids, the notion of outerness in the groupoid situation, to reach  a double crossed product theorem in our context. 

The fourth  chapter deals with match pairs of groupoids, the quantum groupoids associated with, and their relation with the algebraic point of view in [AA].

So a natural prolongement of this article will be the generalization of these constructions, and other extensions of the match pair procedure,  in the direction of Lesieur's locally compact groupoids . An other will be a characterization of these objects in terms of cleft extensions in the spirit of S.Vaes and L.Vainerman [VV].

I want to thank a lot L.Vainerman who suggested me the  actions of groupoids direction,  D.Bisch for some explanations about outer and proper actions of groups, and also M.Enock, S.Baaj and M.C.David for the numerous discussions we had.

\section{ Multiplicative partial isometries and quantum groupoids }

\subsection{ Multiplicative partial isometries}
\label{intro}

\subsubsection{\bf{Notations}}
\label{notation}
In this article,  $N$ is a finite dimensional von Neumann algebra, so 
$N$ is isomorphic to a sum of matrix algebras $\underset\gamma \oplus
M_{n_\gamma}$, we denote  the family of minimal central 
projections of $N$ by $\{ p^\gamma \}$, and we denote  a given family of matrix
units for $N$ by $\{e^\gamma_{i,j} / 1 \leq
i,j
\leq n_\gamma 
\}$. We shall denote  the  opposite  von Neumann algebra of
$N$ by $N^o$, so this is  $N$ with the opposite multiplication, hence a matrix
unit of
$N^o$ is given by the transposed  of $N's$ : $\{e^\gamma_{j,i} / 1 \leq i,j \leq
n_\gamma \}$. The element $f = \underset \gamma \Sigma
\underset {i,j }\Sigma
\frac{1}{n_\gamma}{e^\gamma_{i,j}}^o\otimes e^\gamma_{j,i} $ is the only
projection  of $ N^o\otimes N$ such that, for any $n$ 
in $N$: $f(  n^o \otimes 1) = f(1 \otimes n)$ and   if
$f(1 \otimes n) = 0 $ then $n=0$.

Let $M_1,M_2$ be two  von Neumann algebras. Let $s$ (resp., $r$) be a faithful non degenerate antirepresentation (resp., a representation) from $N$ 
to $M_1$ (resp., $M_2$), then $s$ can be viewed as a representation 
$s^o$ of $N^o$. Let us define:

$$e_{s,r} = (s^o \otimes r)(f) =  \underset \gamma \Sigma \underset
{i,j }\Sigma \frac{1}{n_\gamma}s(e^\gamma_{i,j})\otimes r(e^\gamma_{j,i}).$$

As an  obvious generalization of  Lemma 2.1.2 in   [Val1], 
$e_{s,r}$ is a projection in
$s(N)\otimes
r(N)$, $e_{s,r}$ is \underline {the only}  projection $e$ in $M_1\otimes 
r(N)$ which satisfies the following two conditions:

\hskip 0.4cm a) For every  $m_1$,$m_2$ in $M_1$ and $M_2$ 
respectively, the relation $e(m_1 \otimes 1) = 0 $ implies $m_1 = 0$ and 
the relation $e(1\otimes m_2)=0$ implies $m_2 = 0$,

\hskip 0.4cm b) for every $n$ in $N$: $e(s(n)\otimes1) = e(1 \otimes r(n))$ .

If $Tr_{H_2}$ is the canonical trace of $H_2$ and if $H_2$ is finite dimension,  $(i\otimes Tr_{H_2})(e_{s,r})$ is
positive and invertible in the center of $s(N)$.
\vskip 0.5cm
Let $H$ be a finite dimensional Hilbert space, and let $\alpha$ (resp., 
$\beta, \widehat \beta$) be an injective non degenerate representation  
(resp., two injective non degenerate antirepresentations) of $N$, which 
commute two by two pointwise. We also suppose that $tr \circ \alpha
= tr \circ \beta  = tr \circ \hat \beta$, where $tr$ is the canonical trace on $H$ 
and let us note  $\tau = tr \circ \alpha$. One must keep in mind  that $\beta$ 
and $\alpha$ are a representation and an  antirepresentation of $N^o$.

\subsubsection{\bf {Definition}}\label{MichelEnock1}
{\it  We call a multiplicative partial isometry with the base
\newline
$(N, \alpha, \beta, \hat \beta)$ every  partial isometry $I$ whose   initial 
(resp., final) support is $e_{\widehat \beta,\alpha}$ (resp.,
$e_{\alpha,\beta}$) and such that:

1) $I$ commutes with $\beta(N) \otimes \widehat \beta(N)$ ,

2) For every  $n,n'$ in $N$, one has: $I(\alpha(n) \otimes \beta(n')) =  (
\widehat\beta(n') \otimes \alpha(n))I$,

3) $I$ satisfies the  pentagonal relation: $ I_{12}I_{13}I_{23} =
I_{23}I_{12}.$}

\vskip 0.5cm

By Lemmas 2.3.1, 2.4.2, 2.4.6 in [Val1], one has:

\subsubsection{{\bf Notations and lemma}}
\label{reg}
{\it Let $I$ be  a multiplicative partial isometry with the base $(N,
\alpha, \beta, \hat \beta)$, let's denote the  set 
$
\{(\omega
\otimes i)(I)/
\omega
\ $linear form on $\mathcal{L}(H) \}$ by $S$,
and let's denote the set
$
\{(i \otimes \omega)(I)/ \omega\ $   
linear form on $\mathcal{L}(H)\}$ by $\widehat S $. Then $S$ and $\widehat S $
are non degenerate subalgebras of $\mathcal L(H)$ .}

\subsubsection{\bf{Lemma}} (cf. Lemme 2.6.2 of [Val1]).
\label{regu0}
{\it $S$ and $\hat S$ are sub von Neumann algebras of $\mathcal L(H)$ 
if and only if  $ \{(i \otimes \omega)(\Sigma
I)/ \omega$ is a linear form on
$\mathcal{L}(H) \} =\alpha(N)'$; in this case one says that $I$ is {\bf regular}.}

\subsection{Quantum groupoids}
\label{decadix}
Let us recall the definition of a quantum groupoid (or a weak Hopf 
$C^\star$ -algebra):

\subsubsection{\bf{Definition}}(G.B\"ohm, K.Szlach\'anyi, F.Nill) 
[BoSzNi]\label{Bohm}

{\it A weak Hopf  $C^*$-algebra is a collection $(A, \Gamma, \kappa,
\epsilon)$ where: $A$ is a finite-dimensional $C^*$-algebra (or von Neumann algebra), $\Gamma: 
A \mapsto A \otimes A$ is a generalized coproduct, which means that: 
$(\Gamma\otimes i)\Gamma = (i \otimes \Gamma)\Gamma$, $\kappa$ is an
antipode on $A$, i.e., a linear application  from $A$ to $A$ such 
that $(\kappa \circ *)^2 = i$ (where $*$ is the involution on $A$), 
$\kappa(xy) = \kappa(y)\kappa(x)$ for every  $x,y$ in $A$ with
$(\kappa \otimes \kappa) \Gamma = \varsigma\Gamma \kappa$ (where 
$\varsigma$ is the usual flip on $A\otimes A$). 

We suppose also that 
$(m(\kappa \otimes i)\otimes i)(\Gamma \otimes i)\Gamma(x) = 
(1 \otimes x)\Gamma(1)$ (where $m$ is the multiplication of
tensors, i.e., $m(a \otimes b) = ab$), and that $\epsilon$ is a   
counit, i.e., a positive linear form on $A$ such that
$(\epsilon \otimes i)\Gamma = (i\otimes \epsilon)\Gamma = i$, and 
for every  $x,y$ in $A$: $(\epsilon \otimes\epsilon) ((x \otimes 1) 
\Gamma(1)(1 \otimes y)) = \epsilon(xy)$. }

\subsubsection{\bf{Results}} (cf. [NV1], [NV3],[BoSzNi])
\label{denuit}
{\it If $(A, \Gamma, \kappa, \epsilon)$ is a weak Hopf  $C^*$-algebra,
then the following assertions are true:

1) The sets 
\[ A_t =\{ x\in A/\Gamma(x) = \Gamma(1)(x\otimes 1) = 
(x \otimes 1)\Gamma(1), \}\] 
\[A_s= \{ x \in A / \Gamma(x) = 
\Gamma(1)(1\otimes x) = (1\otimes x)\Gamma(1)\} \] 
are commuting sub $C^*$-algebras of $A$ and $ \kappa(A_t) = A_s$; one calls them  respectively target and 
source  Cartan subalgebra of $(A, \Gamma, \kappa, \epsilon)$.

2) There exists a unique faithful positive linear form $\phi$, called 
the normalized Haar measure of $(A,\Gamma,\kappa,\epsilon)$, satisfying
the following three properties:

$\phi\circ\kappa = \phi$, $(i \otimes \phi)(\Gamma(1)) = 1$ and, for every 
$x,y$ in $A$:
\[(i\otimes\phi)((1\otimes y)\Gamma(x)) = \kappa((i\otimes
\phi)(\Gamma(y)(1\otimes x))).\] 

3) The application $E^s_\phi = (\phi \otimes i)\Gamma$  (resp., $E^t_\phi = 
(i\otimes \phi)\Gamma$) is the conditional expectation 
with values in the source (resp., target) Cartan subalgebra, such that:
$\phi\circ  E^s_\phi  = \phi$ (resp., $\phi \circ  E^t_\phi  = \phi$), it 
is called a  source (resp., target) Haar conditional expectation. If $g_s = E^s_\phi(p)^\frac{1}{2}$ and $g_t = E^t_\phi(p)^\frac{1}{2}$, then one has 
$g_t =\kappa(g_s)$. For every $a$ in
$A$, $ \kappa^2(a) = g_tg_s^{-1}ag_t^{-1}g_s$, and the modular
group $\sigma_{-i}^\phi$ is given by
$\sigma_{-i}^\phi(a) = g_tg_sag_t^{-1}g_s^{-1}$; this leads to a polar
decomposition $\kappa  = j \circ Adg_t$ , where
$j$ is the involutive anti-homomorphism of $A$ (coinvolution) defined by 
$j(y) = g_s\kappa(x)g_s^{-1}$ for any $x$ in $A$..}

\vskip 0.5cm

It is shown in [Val1], that  $I$ is regular, if and only if it generates two quantum 
groupoids in duality:

\subsubsection{\bf{Proposition}}
\label{oct99}
{\it If  $I$ is  regular, then  one can define two quantum groupoids in duality $(S,\Gamma, \kappa,
\epsilon)$ and $(\widehat S, \widehat\Gamma, \widehat \kappa,
\widehat\epsilon)$ , 
by the formulas:

for any $s \in S : \Gamma(s) = I(s \otimes 1)I^*$, 

for any $\hat s \in \hat S : \hat \Gamma(\hat s) = I^*(1\otimes \hat s)I$, 

for any $\omega \in \mathcal L(H))_*$
$\kappa(\omega \otimes i)(I)) = (\omega \otimes i)(I^*)$,
$\hat \kappa((i \otimes \omega)(I) = (i \otimes \omega)(I^*)$,

$\epsilon ((\omega \otimes i)(I)) = \omega(1)$ and
$\hat \epsilon ((i \otimes \omega)(I)) = \omega(1)$. 

These two quantum groupoids are in duality, using the following
bracket: 

$<(\omega\otimes i)(I), (i \otimes \omega')(I)> = (\omega \otimes \omega')(I)$. 

}

\vskip 1.2cm
This  proposition has a reciproque (see [Val] 3.2.3 and [Val3]):
\vskip 0,8cm

\subsubsection{\bf{Proposition}}
\label{oct00}
{\it Let $(A,\Gamma, \kappa,\epsilon)$ be any quantum groupoid such that $\kappa$ is involutive on Cartan subalgebras, and let $\phi $ be it's normalized Haar measure, then with  GNS notations, the application $I:  \Lambda_{\phi\otimes \phi} (x \otimes y) \mapsto  \Lambda_{\phi\otimes \phi}\big(  \Gamma(x)(1 \otimes y)\big)$ is a regular mpi on $H_\phi \otimes H_\phi$, and the GNS representation $\pi_\phi$ is a isomorphism of $(A,\Gamma, \kappa,\epsilon)$ on the quantum groupoid $(S,\Gamma, \kappa,
\epsilon)$ given by proposition \ref{oct99} }

\vskip 0,8cm

In fact, in [Val2], we proved that any regular mpi $I$ can be modified in an irreducible form for which there exist  separating and cyclic vectors $e$ and $\hat e$ for $S$ and $\hat S$ respectively, so that these are in standard form on $H$. With Tomita's theory notations, the unitary  $U = J\hat J = \hat JJ$  leads to define a fourth representation $\hat \alpha = Ad U \circ \alpha: N \to \mathcal L(H)$, and two new mpi $\hat I = \Sigma(U \otimes 1)I(U \otimes 1)\Sigma$, $\tilde I = \Sigma(1 \otimes U)(I(1 \otimes U)\Sigma$  over the base $(N^o,\beta, \hat \alpha, \alpha)$ and $(N^o, \hat \beta, \alpha, \hat  \alpha)$. Applying proposition \ref{oct99} to $\hat I$ (resp. $\tilde I$) leads to define   an other quantum groups on $\hat S'$ (resp. $S'$) the commutant  in $\mathcal L(H)$ of $\hat S$ (resp. $S$), let's denote them $(S',\Gamma',\kappa', \epsilon')$ and $(\hat S', \hat \Gamma',\hat \kappa', \hat \epsilon')$. With these notations one has:

\subsubsection{\bf{Proposition (four corners lemma) }}
\label{corner}
{\it One  has $S\hat S = \hat \alpha(N)' = \{s\hat s / s \in S, \hat s  \in \hat S \}$ (the{\bf  Weyl algebra}),  $S \cap \hat S = \alpha(N)$, $S \cap \hat S' = \beta(N)'$, $ S' \cap \hat S = \hat \beta(N)$ and $S' \cap \hat S' = \hat \alpha(N)$. }

\subsubsection{\bf{Notation }}
\label{weyl}
{\it We shall denote by $\mathfrak W$, the Weyl algebra $S \hat S$. }

\subsection{ The commutative example}
\label{metz}

\vskip 1cm
Let's recall that a groupoid $\mathcal G$ is a small category the  morphisms of which are all invertible. In  all what follows, $\mathcal G$ is finite. One can assimilate the set of object, noted  $\mathcal G^0$, to a subset of the morphisms. So a (finite)  groupoid can also be viewed as a set $\mathcal G$ together with a not everywhere defined  multiplication for which there is a set of unities $\mathcal G^0$, two applications , source denoted by $s$ and range by $r$, from $ \mathcal G$  to $\mathcal G^0$ so that the product $xy$ of two elements $x,y \in \mathcal G $ exists  if and only if $s(x) = r(y)$;  every element $x \in \mathcal G$ has a unique inverse $x^{-1}$, and one has $x(yz) = (xy)z$ whenever the two members have a sense. We refer to [R] for  the fondamental structures and notations for groupoids.

Let's denote $H = l^2(\mathcal G)$, with the usual notations. Actually there exists four natural irreducible mpi associated to $\mathcal G$ ( for an other example
see  [Val1] ¤4.1). Let $I_\mathcal G$ be the  mpi defined  for any $x,y \in \mathcal G$,
$\xi
\in l^2(\mathcal G)$, by:

$ I_\mathcal G(\xi)(x,y) = \xi(xy,y)$ if $s(x) = r(y)$ and $ I_\mathcal G(\xi)(x,y)
= 0$ otherwise. 

Here, $\alpha = \hat \beta$ and $\hat \alpha =  \beta$,  which are given by the source and target 
functions $s$ (resp. $r$), so for every $n \in N$: $\alpha(n) =
\hat \beta(n) = s \circ n$ and $\hat \alpha(n) =  \beta(n) = r
\circ n$. One has  $S= C(\mathcal G)$ and $N = C(\mathcal G^0)$ which are the
commutative involutive algebras of complex valued functions on
$\mathcal G$ and $\mathcal G^0$ respectively, and $\hat S = \mathcal R(\mathcal G) $ = $ \{  \sum_{x \in {\mathcal G}} a_x \rho (x )\} $ (the right regular algebra of $\mathcal G$)  where $\rho(x) $ is the  partial isometry given by the formula $(\rho(x)\xi)(t) = \xi(tx)$ if $x \in \mathcal G^{s(t)}$ and $= 0$ otherwise,  $\hat S' = \mathcal L(\mathcal G)= \{  \sum_{s \in {\mathcal G}} a_s \lambda(s) \} $  ( the left regular algebra of $\mathcal G$), where $\lambda(s) $ is the  partial isometry given by the formula $(\lambda(s)\xi)(t) = \xi(s^{-1}t)$ if $t \in \mathcal G^{r(s)}$ and $= 0$ otherwise, $S' = S$. 
The two  $C^*$-quantum groupoids structures on $S$ and $\hat S$ are given by \vskip 0,5cm
\begin{itemize}

\item \   \    \   \    \    \  {\bf Coproducts}:

$\Gamma_{\mathcal G}(f)(x,y) = f(xy)$ if  x,y  are composables and $f \in C(\mathcal G)$.

 \   \    \   \    \    \   \   \    \   \    \    \  \   \    \   \     \   \    \    \    $=0$   \   \     \   \    \    \   \  \  otherwise

$\hat \Gamma_{\mathcal G}(\rho(s)) =  \rho(s) \otimes \rho(s) $

\item \   \    \   \    \    \  {\bf Antipodes:}

$\kappa_{\mathcal G}(f)(x) = f(x^{-1}),$  \  \  $Ê\hat \kappa_{\mathcal G}(\rho(s)) =  \rho(s^{-1}) =  \rho(s)^* $

\item \   \    \   \    \    \  {\bf Counities:}

$\epsilon_{\mathcal G}(f) = \sum_{u\in \mathcal G^0}f(u)$,   \  \  \  $ \hat \epsilon_{\mathcal G}(\rho(s)) = 1$

\end{itemize}

\section{Actions of  quantum groupoids on  von Neumann modules}

The aim of this section is to give a framework for actions of quantum groupoids, for further extrapolation of  [N]  in the von Neumann algebras context and of [Y] in the quantum groupoids one. Our definitions  have direct generalizations to infinite dimension ( see [EV] definition 7.1) in the von Neumann  algebras context.

In all what follows $I$ will be an irreducible regular mpi over the base $(N,\alpha, \beta, \hat \beta)$, we shall use the  notations of paragraph \ref{decadix}, in particular one has: $S_s =\hat S_t = \alpha(N)$, $S_t = \beta(N)$ and $\hat S_s = \hat \beta(N)$,  and by [Val2] 3.1,  $I$ leads to define two  other mpi $\hat I$ and $ \tilde I $, so one has   two other quantum groupoids for the commutants $ S'$ and  $\hat S'$.....). Let  $A$ be a von Neumann acting on an hilbert space $H$.

\subsection{\bf{Actions of  quantum groupoids}}
\subsubsection{\bf{Notations}}
\label{projo}
{Let $b$ be any unital faithful anti-representation  $S_t   \to A$,  let $i$ be the canonical inclusion $S_t \to S$, we shall denote by $e_{b,i}$ the projection associated to this situation by \ref{notation} . Let $b'$ be any unital faithful representation $A_t   \to A$ and $\kappa$ viewed as a restriction $S_t \to S_s$, we shall denote by  $e_{b',\kappa}$ the projection associated with this by \ref{notation} (applied to $S_t^o$).

\subsubsection{\bf{Definition}}
\label{debut}
{\it With notations above,  let $b$ be a unital faithful anti-representation  (resp.representation) $S_t  (= \beta(N)) \to A$ , hence $A$ appears to be a right (resp. left) module over $S_t$, one calls a {\bf right } (resp.{\bf left}) action of $(S, \Gamma,\kappa)$ on $(A,b)$,  any application $\delta$ (resp.$\gamma)$) such that:

1) $\delta$ (resp. $\gamma$) is an injective normal homomorphism  $ A \to A \otimes S$ (not unital in general),

2) $(\delta \otimes i)\delta = ( i \otimes \Gamma)\delta$ (resp. $(\gamma \otimes i)\gamma= ( i \otimes \varsigma \Gamma)\gamma$),

3)  for any $x \in S_t: \delta(b(x)) = e_{b,i}(1 \otimes \kappa(x))$

(resp.for any $x \in S_t : \gamma(b'(x)) = e_{b', \kappa}(1 \otimes x)$) }

\subsubsection{\bf{Definition}}
\label{codebut}
{\it We shall call right (resp. left) coaction of $(S,\Gamma,\kappa)$, any right (resp.left) action of $(\hat S, \hat \Gamma, \hat \kappa)$.}

\subsubsection{\bf{Remarks}}
\label{attention}
{\it 1) If $S_t \subset Z(S)$ (the center of $S$), then, as for every $a$ in $A$ an $x$ in $S_t$ one has: $\delta(ab(x)) = \delta(b(x)a)$, hence $b(S_t) \subset Z(A)$, in particular, when it is not a  quantum group, such a quantum groupoid   can not act on a factor. 

2) If $\alpha$ is a left action of $(S,\Gamma,\kappa )$ on $(A,b)$ it appears to be a right action of $(S, \varsigma \Gamma, \kappa) $, associated with $\Sigma \hat I^*\Sigma$, on $(A , b \circ \kappa)$.

3) $\Gamma$ is a right action of $(S,\Gamma, \kappa)$ on $(S,\kappa \mid \beta(N))$}

\vskip 0,4cm
Hence till the end we shall only deal with right actions.
\vskip 0,4cm

\subsubsection{\bf{Lemma and definition}}
\label{invariants}
{\it Let $\delta$ be a right  action of $(S, \Gamma, \kappa)$ on $(A,b)$, then one has:
$$  \{ a \in A /  \delta(a) = e_{b,i}(a \otimes 1) \}  =  \{ a \in A \cap b(S_t)' /  \delta(a) = e_{b,i}(a \otimes 1) \},$$ this is a von Neumann subalgebra of $A \cap b(S_t)' $  we shall call the fixed point subalgebra of $\delta$ and note it $A^\delta$.}
\vskip 0.5cm
\begin{dm}
As $S_t$ and $\kappa(S_t)$ commute, hence for every $y$ in $S_t$, one has:

$\delta(b(y)) =  e_{b,i}(1\otimes \kappa(y) ) = (1\otimes \kappa(y))e_{b,i} = e_{b,i}(1\otimes \kappa(y))e_{b,i}$, hence for any $a$ in $A$, such that $\delta(a) = e_{b,i}(a \otimes 1)$, one has: $\delta(ab(y)) = e_{b,i}(a \otimes \kappa(y))e_{b,i} = \delta(b(y)a)$, the lemma follows. \end{dm}

\subsubsection{\bf{Proposition }}
\label{tdelta}
{\it Let $\delta$ be a right  action of $(S, \Gamma, \kappa)$ on $(A,b)$, then the application:
$$T_\delta = ( i \otimes \phi)\delta$$ 
is a faithful conditional expectation $A \to A^\delta$  .}
\vskip 0.5cm
\begin{dm} If $(n_j)$ is a matrix unit for $N$, then $(\beta(n_j))$ (resp. $\alpha(n_j)$) is also a matrix unit for $S_t$ (resp. $\hat S_t$), so, for any $a \in A $,  by \ref{notation} one has:
\begin{align*}
T_\delta(a)
&= (i \otimes \phi)(e_{ b,  i}(a \otimes 1) )
= (i \otimes \phi )(\underset j \sum b(\beta (n_j)) \otimes \beta(n^*_j))a  \\
&= ( b \circ \kappa)(i \otimes \phi)(\Gamma(1) )a =  ( b \circ \kappa)(1)a = a
\end{align*}

For every $a \in A$, let's use the notations $\delta(a) = a_1 \otimes a_2$, one has:
\begin{align*}
\delta(T_\delta(a))
&= \delta((i \otimes \phi)\delta(a)) = (i \otimes i \otimes \phi)(\delta \otimes \delta)\delta(a) = (i \otimes i \otimes \phi)(i \otimes \Gamma)\delta(a) \\
&= (i \otimes ( i \otimes \phi)\Gamma)\delta(a) = (i \otimes E_t)\delta(a) = (i \otimes E_t)(\delta(1)\delta(a)) = e_{b,i}(i \otimes E_t)\delta(a) \\
&= e_{b,i}(a_1 \otimes E_t(a_2)) = e_{b,i}(a_1b(E_t(a_2)) \otimes 1)
\end{align*}

So one deduces that:
\begin{align*}
(i \otimes \phi)\delta(T_\delta(a)) 
&= (i \otimes \phi)(e_{b,i}(a_1b(E_t(a_2)) \otimes 1)) = (i \otimes \phi)(e_{b,i})a_1b(E_t(a_2))\\
&= a_1b(E_t(a_2))
\end{align*}

but in an other hand:
\begin{align*}
(i \otimes \phi)\delta(T_\delta(a)) 
&= (i \otimes \phi)(e_{b,i}(a_1 \otimes E_t(a_2))) = (i \otimes \phi)(\delta(1)\delta(a)) = (i \otimes \phi)\delta(a) = T_\delta(a) 
\end{align*}
Hence one has: $T_\delta(a) =  a_1b(E_t(a_2))$, replacing this in the expression of $\delta(T\delta(a))$, one has: 
$$ \delta(T_\delta(a)) =  e_{b,i}(a_1b(E_t(a_2)) \otimes 1) = e_{b,i}(T_\delta(a) \otimes 1)$$
this implies that $T_\delta(A)  =  A^\delta$.

Using the fact that any element of $A^\delta$ commutes with $b(S_t)$, for any $b,c$ in $A^\delta$ and any $a$ in $A$, one has:
\begin{align*}
T_\delta(cab)
&= ( i \otimes \phi)\delta(cab) = ( i \otimes \phi)(\delta(c)\delta(a)\delta(b)) =( i \otimes \phi)(  e_{b,i}(c \otimes 1)\delta(a)e_{b,i}(b \otimes 1)) \\
&= ( i \otimes \phi)(  (c \otimes 1)e_{b,i}\delta(a)e_{b,i}(b \otimes 1)) = c( i \otimes \phi)(  \delta(1)\delta(a)\delta(1) )b \\
&= cT_\delta(a)b
\end{align*}
The proposition follows.
\end{dm}

\subsubsection{\bf{Definition}}
\label{debut1}
{\it Let $\delta$  (resp.$\hat \delta$) be a right  action (resp.coaction) of  $(S, \Gamma, \kappa)$ on a von Neumann module  $(A,b)$ (resp.$(A,\hat b)$, the crossed product $ A \underset \delta \rtimes S$ (resp . $ A \underset {\hat \delta} {\hat  \rtimes} \hat S$)  the sub-von Neumann algebra of $e_{b,i}(A \otimes \mathcal L(\mathcal H))e_{b,i}$ (resp. $e_{\hat b,\hat i}(A \otimes \mathcal L(\mathcal H))e_{\hat b, \hat i}$) generated by $\delta(A)$ and $e_{b,i}(1 \otimes \hat S)$ (resp $\hat \delta(A)$ and $e_{\hat b,\hat i}(1 \otimes   S'))$. }

\subsubsection{\bf{Remarks}}
{\it 1) One must keep in mind that the crossed product is degenerated in $A \otimes \mathcal L(\mathcal H)$, and it's unit element is $e_{b,i}$.

2)  If $\alpha$ is a left action, as it is a right action  one can define also a crossed product.

3) As a matter of facts, in this Baaj and Skandalis formalism, $\hat S$ is not equal to the one given by Vaes' theory, but our  crossed product do generalize the  quantum groups one. }

\subsubsection{\bf{Lemma}}
\label{reconquete}
{\it The crossed product $ A \underset \delta \rtimes S$  is the sub  vector space of  $e_{b,i}(A \otimes \mathcal L(\mathcal H))e_{b,i}$ (resp. $e_{b,\kappa}(A \otimes \mathcal L(\mathcal H))e_{b,\kappa}$ generated by the products $\delta(a)(1 \otimes \hat b)$ $a \in A$, $\hat b \in \hat S$ .}
\vskip 0.5cm
\begin{dm}
The equality $(\delta \otimes i)\delta = (i \otimes \Gamma)\delta$, leads to the fact that for any $a \in A$, one has: $ I_{23}(\alpha(a) \otimes 1)I_{23}^* \in \delta(A)\otimes S$ , so as $I_{23}^*I_{23} = e_{\hat \beta, \alpha}$, there exist $a_1, ...,a_k $ in $A$ and $b_1,...,b_k$ in $S$ such that:
\begin{align*}
I_{23}(\delta(a) \otimes 1)= \underset i \sum ( \delta(a_i) \otimes c_i)I_{23}
\end{align*}

Hence for any linear form $\omega$ on $\mathcal L(\mathcal H)$, one has:
\begin{align*}
(1\otimes (i \otimes \omega)(I))\delta(a)
&=  (i \otimes i \otimes \omega)(I_{23}(\delta(a)\otimes1))  =  (i \otimes i \otimes \omega)(\underset i \sum ( \delta(a_i) \otimes c_i)I_{23}) \\
&= \underset i \sum  \delta(a_i) (1 \otimes (i \otimes (\omega.c_i))(I))
\end{align*}
This leads to the lemma .
\end{dm}

\subsubsection{\bf{Notation}}
\label{recrut}
Let's denote $I' = \Sigma (\hat{\hat I})^* \Sigma = \Sigma (\tilde{\tilde I})^* \Sigma = \Sigma(U \otimes U)I^*(U \otimes U)\Sigma$, hence $I'$ is a mpi belonging to $S' \otimes \hat S'$ over the base $(N,\hat \alpha, \beta, \hat \beta)$, which gives on $\hat S'$ the opposite  coproduct $ \hat \Gamma'^{opp} = \varsigma \hat \Gamma'$.

\subsubsection{\bf{Proposition}}
\label{reconquete2}
{\it i) Let $\delta$  be a right action  of $(S,\Gamma,\kappa)$ on$(A,b)$, let $ \hat b$ be the application $\hat S_t (= \alpha(N)) \to  A \underset \delta \rtimes S$ defined for every $n \in N$ by:
$$ \hat b(\alpha(n)) =   e_{b,i}(1 \otimes \hat \beta(n))$$
let  $\hat \delta$  be the application defined for every $x \in  A \underset \delta \rtimes S$ by:

$$\hat \delta(x) =   {\tilde I}_{23}(x \otimes 1) {\tilde I}_{23}^*$$
then $(\hat \delta, \hat b)$ is a (right) coaction of $ (S,\Gamma,\kappa)$ on $( A \underset \delta \rtimes S, \hat b)$ and  $ (A \underset \delta \rtimes S)  \underset {\hat \delta} {\hat \rtimes} \hat S$ is isomorphic to the sub von Neumann algebra of $\delta(1)(A \otimes \mathcal L(\mathcal H))\delta(1)$ generated by $\delta(A)$ and $\delta(1)(1 \otimes \hat SS')$$ (= \delta(1)(1 \otimes \beta(N)')$. }
\vskip 0.5cm
\begin{dm}
By [Val3] 3.1,  the initial support of $ \tilde I$ is $e_{\hat \alpha , \hat \beta }$, but $\hat \alpha(N)$ commutes with $S$ and $\hat S$, this implies that $\hat \delta$ is a normal  homomorphism on $A \underset \delta \rtimes S$; as $ \tilde I \in  S' \otimes \hat S$ and his final support is $e_{\hat \beta,\alpha}$,  for any $a$ in $A$ one has: 
$$\hat \delta(\delta(a))  = (1 \otimes e_{\hat \beta,\alpha})(\delta(a) \otimes 1) $$
 and due to [Val2] proposition 3.1.4, for every $\hat b$ in $\hat S$, one has: 
 $$\hat \delta(e_{b,i}(1\otimes \hat b)) = (e_{b,i} \otimes 1)(1 \otimes \hat \Gamma(\hat b)).$$ So $\hat \delta$ takes values in $A \underset \delta \rtimes S \otimes \hat S$, and for every $x \in \hat S_t$, one has:
\begin{align*}
\hat \delta(\hat b(x))
&= \hat \delta(e_{b,i}(1 \otimes \hat  \kappa(x))) = (e_{b,i} \otimes 1)(1 \otimes \hat \Gamma(\hat \kappa(x))) \\
&= (e_{b,i} \otimes 1)(1 \otimes e_{\hat \beta,\alpha})(1 \otimes 1 \otimes \hat \kappa(x))
\end{align*}

But , if $(n_j)$ is a matrix unit for $N$, then $(\beta(n_j))$ (resp. $\alpha(n_j)$) is also a matrix unit for $S_t$ (resp. $\hat S_t$), so by \ref{notation} one has:
\begin{align*}
e_{\hat b, \hat i} 
&= \underset j \sum \hat b(\alpha (n_j)) \otimes \alpha(n^*_j)  = \underset j \sum e_{b,i}(1  \otimes \hat  \beta (n_j^*) \otimes \alpha(n_j)) \\
&= (e_{b,i} \otimes i) (1  \otimes \underset j \sum \hat \beta(n_j^*) \otimes \alpha(n_j) ) \\
&=  (e_{b,i} \otimes 1)(1 \otimes e_{\hat \beta,\alpha}),
\end{align*}
replacing this in the previous calculus:
\begin{align*}
\hat \delta(\hat b(x))
&= e_{\hat b, \hat i} (1_{A \underset \delta \rtimes S} \otimes \hat \kappa(x))
\end{align*}
 Now let's verify condition 2) of definition \ref{debut}
 \begin{align*}
 (\hat \delta \otimes i)\hat \delta(\delta(a))
 &= (\hat \delta \otimes i)((1 \otimes e_{\hat \beta,\alpha})(\delta(a) \otimes 1)) = (\hat \delta \otimes i)((1 \otimes e_{\hat \beta,\alpha})(\delta(1) \otimes 1)(\delta(a) \otimes 1)) \\
 &= \underset j \sum (\hat \delta \otimes i)((1 \otimes \hat \beta(n_j) \otimes \alpha(n_j^*))(e_{b,i} \otimes 1)(\delta(a) \otimes 1))  \\
 &= \underset j \sum (\hat \delta \otimes i)((e_{b,i} (1 \otimes \hat \beta(n_j)) \otimes \alpha(n_j^*))(\delta(a) \otimes 1))   \\
 &= \underset j \sum (\hat \delta \otimes i)(e_{b,i} (1 \otimes \hat \beta(n_j) \otimes \alpha(n_j^*)))(\hat \delta \otimes i)(\delta(a) \otimes 1)   \\
 &= \underset j \sum \hat \delta (e_{b,i} (1 \otimes \hat \beta(n_j) )\otimes \alpha(n_j^*))(\hat \delta (\delta(a) )\otimes 1) \\
 &= \underset j \sum ((e_{b,i} \otimes 1)(1 \otimes \hat \Gamma( \hat \beta(n_j)))\otimes \alpha(n_j^*))((1 \otimes e_{\hat \beta,\alpha})(\delta(a) \otimes 1)\otimes 1)  \\
 &= \underset j \sum ((e_{b,i} \otimes 1)(1 \otimes e_{\hat \beta,\alpha}(1 \otimes \hat \beta(n_j)))\otimes \alpha(n_j^*))((1 \otimes e_{\hat \beta,\alpha})(\delta(a) \otimes 1)\otimes 1)  \\
 &= \underset j \sum ((e_{b,i} \otimes 1)(1 \otimes e_{\hat \beta,\alpha}(1 \otimes \hat \beta(n_j)))\otimes \alpha(n_j^*))((1 \otimes e_{\hat \beta,\alpha})(\delta(a) \otimes 1)\otimes 1) \\
 &= (e_{b,i} \otimes 1 \otimes 1)(1 \otimes e_{\hat \beta,\alpha} \otimes 1)(1 \otimes 1 \otimes  e_{\hat \beta,\alpha} )((1 \otimes e_{\hat \beta,\alpha})(\delta(a) \otimes 1)\otimes 1)   \\
 &= (e_{b,i} \otimes 1 \otimes 1)(1 \otimes e_{\hat \beta,\alpha} \otimes 1)(1 \otimes 1 \otimes  e_{\hat \beta,\alpha} )(\delta(a) \otimes 1 \otimes 1) 
 \end{align*}
 and on the other side:

 \begin{align*}
 (i \otimes \hat \Gamma)\hat \delta(\delta(a))
 &=  (i \otimes \hat \Gamma)((1 \otimes e_{\hat \beta,\alpha})(\delta(a) \otimes 1)) = (\hat \delta \otimes i)((1 \otimes e_{\hat \beta,\alpha})(\delta(1) \otimes 1)(\delta(a) \otimes 1)) \\
 &= \underset j \sum (i \otimes \hat \Gamma)((1 \otimes \hat \beta(n_j) \otimes \alpha(n_j^*))(e_{b,i} \otimes 1)(\delta(a) \otimes 1))  \\
 &= \underset j \sum (i \otimes \hat \Gamma)((e_{b,i} (1 \otimes \hat \beta(n_j)) \otimes \alpha(n_j^*))(\delta(a) \otimes 1))   \\
 &= \underset j \sum (i \otimes \hat \Gamma)(e_{b,i} (1 \otimes \hat \beta(n_j) \otimes \alpha(n_j^*)))(i \otimes \hat \Gamma)(\delta(a) \otimes 1)   \\
 &= \underset j \sum (e_{b,i} (1 \otimes \hat \beta(n_j) )\otimes \hat \Gamma(\alpha(n_j^*)))(\delta(a) \otimes e_{\hat \beta,\alpha}) \\
 &= \underset j \sum (e_{b,i} (1 \otimes \hat \beta(n_j) \otimes e_{\hat \beta,\alpha}(\alpha(n_j^*)\otimes 1)))(\delta(a) \otimes e_{\hat \beta,\alpha}) \\
 &= (e_{b,i} \otimes 1 \otimes 1)(1 \otimes e_{\hat \beta,\alpha} \otimes 1)(1 \otimes 1 \otimes  e_{\hat \beta,\alpha} )(\delta(a) \otimes 1 \otimes 1) \\
 &= (\hat \delta \otimes i)\hat \delta(\delta(a))
  \end{align*}

Also for the others generators, and using Sweedler notations, one has:

\begin{align*}
 (\hat \delta \otimes i)\hat \delta(e_{b,i}(1\otimes \hat b)) 
 &= (\hat \delta \otimes i)((e_{b,i} \otimes 1)(1 \otimes \hat \Gamma(\hat b))) = (\hat \delta \otimes i)((e_{b,i} \otimes 1)(1 \otimes\hat b_1\otimes \hat b_2)) \\
 &=(\hat \delta \otimes i)(e_{b,i} (1 \otimes\hat b_1)\otimes \hat b_2)
 = (e_{b,i} \otimes 1)(1  \otimes \hat \Gamma(\hat b_1))\otimes \hat b_2  \\
 & = (e_{b,i} \otimes 1\otimes 1)(1 \otimes (  \hat \Gamma \otimes i)\hat \Gamma(\hat b))  = (e_{b,i} \otimes 1\otimes 1)(1 \otimes ( i\otimes \hat \Gamma)\hat \Gamma(\hat b))  \\
  &= (e_{b,i} \otimes 1)(1  \otimes \hat b_1 \otimes\hat \Gamma( \hat b_2))   = (i \otimes \hat \Gamma)((e_{b,i} \otimes 1)(1  \otimes \hat \Gamma( \hat b)))  \\
 &= (i \otimes \hat \Gamma)\hat \delta(e_{b,i}(1\otimes \hat b)) 
 \end{align*}
 One can deduce that: $ (\hat \delta \otimes i)\hat \delta = (i \otimes \hat  \Gamma)\hat \delta$.

Now let define the application the one to one morphism $\gamma$ defined on $A \otimes  \mathcal L( \mathcal H)$ by:
$$  \gamma(x) = I^*_{23}(\delta \otimes i)(x)I_{23},  \  \  \  \   \forall x \in A \otimes  \mathcal L( \mathcal H) $$

Obvious calculations give that, for any $a$ in $A$, $\hat s $ in $\hat S$, $d$ in $S'$, one has:

$\gamma(\delta(a)) = \hat \delta(1)(\delta(a) \otimes 1)$, $\gamma(\delta(1)(1\otimes \hat s ) = \hat \delta(1)(1 \otimes \hat \Gamma(\hat s))   $, $\gamma(\delta(1)(1\otimes  s ) = \hat \delta(1)(1 \otimes 1 \otimes s')  $

Hence, $\gamma$ is an  isomorphism  between  the sub von Neumann algebra of $\delta(1)(A \otimes \mathcal L(\mathcal H))\delta(1)$ generated by $\delta(A)$ and $\delta(1)(1 \otimes \hat SS')$$ (= \delta(1)(1 \otimes \beta(N)')$, and $ (A \underset \delta \rtimes S)  \underset {\hat \delta} {\hat \rtimes} \hat S$.

\end{dm}

\subsection{Actions of groupoids}
Let's explain what is  an action $\alpha$ of the commutative quantum groupoid $(C(\mathcal G), \Gamma_{\mathcal G}, \kappa_{\mathcal G})$ where  $\mathcal G$ is any finite groupoid, on  a von Neumann module $(A,b)$. In fact, the application $b$ is clearly equivalent to the given of a decomposition  $A = \underset {u \in \mathcal G^0} \oplus A_u$, where each $A_u$ is a von Neumann algebra, the relation is given, for every $u \in \mathcal G^0$, by  $b(\delta_u) = 1_u$, where $\delta_u$ is the Dirac fonction for $u$ and $1_u$ the identity element of $A_u$ (a projection in $Z(A)$). Hence $A$ appears to be a module over $C(\mathcal G^0)$. 

\subsubsection{\bf Definition}
\label{action}
{\it An action of $\mathcal G$ on $A$ is any covariant  functor from the category $\mathcal G$ to the category whose objects are the element of the set  $\{ A_u ,  \   u \in \mathcal G^0 \}$ and the morphisms the von Neumann algebras isomorphisms.}

Hence, for any $g \in \mathcal G$, it exists a morphism $\alpha_g:  A_{s(g)} \mapsto A_{r(g)} $,  in order that for any pair $(g,g')$ of composable elements, one has: $\alpha_{gg'} = \alpha_g \alpha_{g'}$

As it can be decomposed  in its connected  classes, we can suppose that $\mathcal G =  \underset i  \bigsqcup \hskip 0.2cm  X_i \times X_i \times G_i$, where $X_i$ is a finite set and $G_i$ is a finite group. In fact one has $\mathcal G^0 =  \underset i  \bigsqcup \hskip 0.2cm X_i$ and $G_i$ is isomorphic to the isotropy group $G_u^u$ for any $u \in X_i$. 

\vskip 0.5cm

\subsubsection{\bf Proposition} 
\label{bete}
{\it Any finite groupoid $\mathcal G$ acts on $(R^{\mathcal G^0},b_{\mathcal G^0})$, where $R$ is the hyperfinite type $II_1$ factor and $b_{\mathcal G^0}: f \mapsto (f(u)1)_{u \in \mathcal G^0}$.}
\vskip 0.5cm
\begin{dm}
Due to the previous remark,  one can suppose that $\mathcal G = X\times X \times G$ where $X$ is a finite set and $G$ is a group. As it is well known, there exists an action $\beta$ (even outer) of $G$ on $R$. Of course $R^{\mathcal G^0}$ can be decomposed in it's cartesian components, each of them is  in fact $R$ itself; up to this identification, one can define for any $(x,y,g) \in \mathcal G$ $\alpha_{(x,y,g)} = \beta_g$, one easily sees that this is an action.
\end{dm}

\vskip 0.5cm

\subsubsection{\bf {Lemma}}
\label{michel}
{\it  For any action $\alpha$ of $(C(\mathcal G), \Gamma_{\mathcal G}, \kappa_{\mathcal G})$ (the image of which can be viewed in $A \otimes \mathcal L(l^2(\mathcal G))$ )  on $(A,b)$, and  for any $g \in \mathcal G$,  one has: $ (1 \otimes  \lambda(t )^*)   \alpha(A) (1 \otimes  \lambda(t)) \subset \alpha(A)$.}
\vskip 0.5cm
\begin{dm} This is just a generalization of the demonstration of Proposition 1.3 ii) in [E3], replacing the unitary $W_G$ by  the adjoint of the regular mpi defined in [Val1] 4.1.
\end{dm}
\vskip 0.5cm
In fact the two notions of  action are equivalent: 

\subsubsection{\bf {Proposition}}
{\it  i) For any action of $\mathcal G$ on $(A,b)$, then the application $\delta_\alpha$ (resp. $ \gamma_\alpha$ ): $A \mapsto A \otimes C(G) ( = C(G, A))$ defined for every $a \in A$ by $\delta_\alpha(a) : g \mapsto \alpha_g(a_{s(g)})$ (resp. $\gamma_\alpha(a) = \alpha_{g^{-1}}(a_{r(g)})$) is a right (resp.left) action of $(C(\mathcal G), \Gamma_{\mathcal G},\kappa_{\mathcal G})$ on $(A,b)$.

ii) For any left (resp. right) action $\gamma$ (resp. $\delta$) of $(C(\mathcal G), \Gamma_{\mathcal G}, \kappa_{\mathcal G})$ on $(A,b)$ there exists a unique action of $\mathcal G$ on $(A,b)$, such that $\gamma= \gamma_\alpha$ (resp. $\delta = \delta_\alpha$).}

\begin{dm} Using lemma \ref{michel} , one can exactly use the arguments in [E3] Proposition 1.3
 \end{dm}
\vskip 0.5cm
A third and more synthetic way to define an action on $(A,b)$ is just to consider the groupoid $Aut (A,b)$ whose base is $\mathcal G^0$ viewed as $ \{ Id A_u, u \in G^0 \}$ and morphisms the isomorphisms $A_u \mapsto A_v$. An action is just a full (i.e. with the same base) subgroupoid of $Aut(A,b)$.
\vskip 0.5cm
\subsubsection{\bf{Remarks}} 
   
1) Our definition of an action agrees with the algebraic definition due to Vainerman and Nikshych in [NV2]: let $\gamma$ (resp. $\delta$) be a left (resp. right) action of $(C(\mathcal G), \Gamma_{\mathcal G}, \kappa_{\mathcal G}, \epsilon_{\mathcal G})$ on $(A,b)$, then if for any $a$ in $A$ and $h$ in $\mathcal G$ one defines $ \lambda(h) \triangleright a = \gamma(a  _{h}$ (resp. $a \triangleleft  \rho(h) = \delta(a)_{h}$), this  is a left (resp.right) action  and the same formula can be used for the inverse assertion.

2) As $\delta(1) = e_{b,r} \not = 1$, $\delta(A)$  is degenerated in $C(G, A)$, so it's better convenient to restrict $\delta$ to $e_{b,r} (A \otimes C(\mathcal G))e_{b,r} $, which can be identified with $\underset {g \in \mathcal G} \oplus A_{r(g)}$, and $\delta(a)$ can be identified with $\big(\alpha_g(a_{s(g)})\big)_g$, in that way, $\delta$ appears to be unital.

\subsubsection{\bf {Notations}}
 Let's consider $H= L^2(A)$, the standard hilbert space of $A$, then $H$ has an orthogonal  decomposition $H = \underset {u \in \mathcal G^0} \oplus H_u$; for $j \in (s,r)$, $e_{b,j} (A \otimes C(\mathcal G))e_{b,j}$ can also be represented as a "diagonal" von Neumann algebra acting on $\underset {g \in \mathcal G} \oplus H_{j(g)}$. For any $h \in \mathcal G$, one can define the operator $(1 _b \otimes_r \rho(h))$ on $\underset {u \in \mathcal G^0} \oplus H_{r(g)}$, by the formula $(1 _b \otimes_r \rho(h))\big ((\xi_g)_{g \in \mathcal G}\big) = (\eta_g)_{g \in \mathcal G} $, where $\eta_g= 0 $ if $s(g) \not = r(g)$ and $\eta_g = \xi_{gh}$ otherwise. For any $a \in A$, one also can define the  operator  $(a _b \otimes_s 1) = \oplus a_{s(g)}$ which acts on $\oplus H_{s(g)}$. Let's also denote by $u_g$ the canonical implementation of $\alpha_g$ for any $g \in \mathcal G$ (theorem 2.18 of [H]), so the operator $U =\oplus u_g$ appears to be a unitary $\oplus H_{s(g)} \to  \oplus H_{r(g)}$. Hence, obviously one has:

 \subsubsection{\bf {Proposition}}
\label{implementation}
{\it The unitary $U$ implements the action $\delta$, for any $a$ in $A$, one has:
$$ \delta(a) = U(a _b \otimes_s 1)U^*.   $$ }

\vskip 1cm
\subsection{Crossed product by groupoids actions and Jones tower}
In all what follows $\alpha$ is an action of $\mathcal G$ on $(A,b)$ and $\delta = \delta_\alpha$, the right action of $(C(\mathcal G), \Gamma_{\mathcal G})$ on $(A,b)$ associated with. So one can consider the inclusion $A^\delta \subset A$. Let's recall Jones basic construction: if $M_0 \subset M_1 $ is an inclusion of von Neumann algebras, and $J$ is the canonical antilinear involutive isometry of $L^2(M_1)$, then one can extend the inclusion by:  $M_0 \subset M_1 \subset M_2 (= JM_0'J)$, that is the basic construction. This paragraph proves simply that $M_2$ is a quotient of the crossed product. First let's give a simple description of this crossed product.

\subsubsection{\bf {Remark}}
{\it  The crossed product of $(A,b)$ by $\mathcal G$ is the sub-von Neumann algebra of $\mathcal L(\underset {g \in \mathcal G} \oplus H_{r(g)})$ generated by  $\delta(A)$ and the operators $(1 _b \otimes_r \rho(h))$.}

\subsubsection{\bf {Lemma}}
\label{orsay}
{\it i) For any $a$ in $A$ and $h$ in $\mathcal G$, one has:

 $(1 _b \otimes_r \rho(h)) \delta(a) =  \delta(\alpha_h(a_{s(h)}))(1 _b \otimes_r \rho(h))$.

 ii) The crossed product $A \underset \delta \rtimes C(\mathcal G)$ is the vector space generated by the products $\delta(a)(1 _b \otimes_r \rho(h))$ for any $(a,h)$ in $A_b \times_r \mathcal G (= \{(a,h) \in A\times  \mathcal G, a \in A_{r(h)} \}$).

iii) $A \underset \delta \rtimes C(\mathcal G)$ is the set of elements in $\mathcal L(\underset {g \in \mathcal G} \oplus H_{r(g)})$, which can be decomposed in a sum of the form $\sum_{h \in \mathcal G} \delta(x^h)(1 _b \otimes_r \rho(h))$, where $x^h \in A_{r(h)}$ for all $h \in \mathcal G$, and this decomposition is unique.}
\vskip 1cm
\begin{dm} For any $a$ in $A$, $h$ in $\mathcal G$ and $(\xi_g)_{g \in \mathcal G}$ in $\underset {g \in \mathcal G} \oplus H_{r(g)}$, one has:

\begin{align*}
(1 _b \otimes_r \rho(h)) \delta(a)\big ( (\xi_g)_{g \in \mathcal G} \big )
&= (1 _b \otimes_r \rho(h))  \big( (\alpha_g(a_{s(g)})\xi_g)_{g \in \mathcal G} \big ) \\
&= (1 _b \otimes_r \rho(h))  ( (\alpha_g(a_{s(g)})\xi_g)_{g \in \mathcal G} \big ) =  (\eta_g)_{g \in \mathcal G} 
\end{align*}
where one has : $  \eta_g = \left \{ \begin{array}{rl}
& \alpha_g(\alpha_h(a_{s(h)}))\xi_{gh} \hskip 0.5cm  \mathrm{if} \   s(g) = r(h) \\
&0 \hskip 0.5cm  \ \mathrm{otherwise}
\end{array}
\right.
$

On the other side For any $b$ in $A$, $h$ in $\mathcal G$ and $(\xi_g)_{g \in \mathcal G}$ in $\underset {g \in \mathcal G} \oplus H_{r(g)}$, one has:
\begin{align*}
 \delta(b)(1 _b \otimes_r \rho(h))\big ( (\xi_g)_{g \in \mathcal G} \big )
&= (\delta(b)\big(\eta'_g)_{g \in \mathcal G}
\end{align*}
where one has : $  \eta'_g = \left \{ \begin{array}{rl}
& \xi_{gh} \hskip 0.5cm  \mathrm{if} \   s(g) = r(h) \\
&0 \hskip 0.5cm  \ \mathrm{otherwise}
\end{array}
\right.
$

hence $ \delta(b)(1 _b \otimes_r \rho(h))\big ( (\xi_g)_{g \in \mathcal G} \big ) = (\alpha_g(b_{s(g)})\eta'_g)_{g \in \mathcal G} = (\eta''_g)_{g \in \mathcal G}$, so one deduces that:
$$  \eta''_g = \left \{ \begin{array}{rl}
&  \alpha_g(b_{r(h)})\xi_{gh} \hskip 0.5cm  \mathrm{if} \   s(g) = r(h)  \\
&0 \hskip 0.5cm  \ \mathrm{otherwise}
\end{array}
\right.
$$

One deduces that $(1 _b \otimes_r \rho(h)) \delta(a) =  \delta(b)(1 _b \otimes_r \rho(h))$ for any $a, b ,h$ such that $b_{r(h)} = \alpha_h(a_{s(h)})$, i) and ii) follow immediatly. 

If one chooses for any $u \in \mathcal G^0$ a base $(a_u^k)$ of $A_u$, one easily sees that the family $(\delta(a_{r(h)}^k)(1 _b \otimes_r \rho(h)))$ is free, hence iii) is a consequence of ii).
\end{dm}

\subsubsection{\bf{Corollary}}
\label{NV}
{\it The application $\sum_{h \in \mathcal G} \delta(x^h)(1 _b \otimes_r \rho(h)) \mapsto x^h \otimes \rho(h)$ leads to an isomorphism of $A \underset \delta \rtimes C( \mathcal G)$ and the corresponding crossed product by L.Vainerman and D.Nikshych.    }

\subsection{Outer actions of groupoids}

\subsubsection{\bf {Definition}}
\label{isotropic}
{\it Let's call isotropic subgroupoid of $\mathcal G$, the subgroupoid of $\mathcal G$, denoted iso($\mathcal G$), equal to  $\{ h \in \mathcal G, s(h) = r(h) \}$. }

\subsubsection{\bf {Remark}}
Obviously iso($\mathcal G$) is the disjoint union of the isotropic groups $\mathcal G^u_u$.
\vskip 1cm

\subsubsection{\bf {Lemma}}
\label{explicit}
{ \it  Let's suppose that $Z(A)$ is isomorphic to $C(G^0)$ (or equivalently each $A_u$ is a factor).  An element $x\in A \underset \delta \rtimes C(\mathcal G)$ commutes with  $\delta(A)$ if and only if for any $h \notin  iso(\mathcal G)$, one has $x^h = 0$ and for any  $h \in iso(\mathcal G)$ with $x^h \not = 0$, $\alpha_h$ is inner with $\alpha_h(a) = (x^h)^{-1}ax^h$, for all $a \in A_{r(h)}$ (hence $x^h$ is invertible).}
\vskip 1cm
\begin{dm}
Due to lemma \ref{orsay}, for any $x \in  A \underset \delta \rtimes C(\mathcal G)$, if $\delta(a)x = x\delta(a)$ for all $a$ in $A$, it means that, for any $a$ in $A$ and $h$ in $G$, one has: $ \delta(a_{r(h)}x^h)(1 _b \otimes_r \rho(h)) = \delta(x^h\alpha_h(a_{s(h)})(1 _b \otimes_r \rho(h))$. Hence, one has:  $a_{r(h)}x^h =  x^h\alpha_h(a_{s(h)})$.

If $h \notin iso(\mathcal G)$, let $a$ be the element $\alpha_{h^{-1}}(x^h)^*$ and so $a_{r(h)} = 0$, then one deduces that: $0 = x^h\alpha_h(a_{s(h)}) = x^h(x^h)^*$, so does $x^h$.

If $h \in iso(\mathcal G)$, then for any $a \in A_{r(h)}$, one has: $ax^h = x^h\alpha_h(a)$, as $A_u$ is  a factor,, one can suppose that $A$ has no trivial weakly closed two side ideal, but $x^hA$ is a two side weakly closed  ideal, so $x^h$ is invertible or equal to zero; the lemma follows.
\end{dm}
\vskip 1cm

\subsubsection{\bf{Remark}} The von Neumann algebra $\delta(A)' \cap A \underset \delta \rtimes C(\mathcal G)$, the relative commutant of $\delta(A)$ in $A \underset \delta \rtimes \mathcal G$, contains $Z(A) \underset \delta \rtimes \mathcal \beta(C(G^0))$, whose elements are of the form:  $\sum_{u \in \mathcal G^0} \delta(x^u)(1 _b \otimes_r \rho(u))$, for $x^u \in Z(A_u)$ which, in the case when $Z(A) = b(C(\mathcal G^0))$,  is just $1_b \otimes_r  \beta(C(G^0))$.
\vskip 1cm

\subsubsection{\bf{Definition}}
{\it i) A  right  action $\delta$ of $C(\mathcal G)$ on a von Neumann module $(A,b)$ is said to be outer if and only if $\delta(A)' \cap A \underset \delta \rtimes C(\mathcal G)$ is equal to $Z(A) \underset \delta \rtimes \mathcal \beta(C(G^0))$.

ii) An  action $\alpha $ of  $\mathcal G$ on a von Neumann module $(A,b)$ is said to be outer  if and only if for any  $h \in iso(\mathcal G)$ such that $h \notin \mathcal G^0$, one has $\alpha_h \in Out A_{r(h)}$.}

\subsubsection{\bf{Remark}} 
{\it The  transitive groupoid $\mathcal G = X \times X \times Out R$ acts outerly on $R^X$.}

\vskip 1cm
As a  consequence of lemme \ref{explicit}, one has:

\vskip 1cm

\subsubsection{\bf{Proposition}}
\label{outerminimal}
{\it In the case when $Z(A) = b(C(\mathcal G^0))$,  any right  action $\delta$ of $C(\mathcal G)$ on a von Neumann module $(A,b)$ is   outer if and only if the action of $\mathcal G$ on $(A,b)$, canonically associated with $\delta$ is outer. }
\vskip 1cm
And finally:

\subsubsection{\bf{Proposition}}
\label{universall}
{\it Any finite groupoid  $\mathcal G$ acts outerly on the von Neumann module $(R^{\mathcal G^0},b_{\mathcal G^0})$. }
\vskip 1cm

\begin{dm}
As in proposition \ref{bete}, one can suppose $\mathcal G =  X\times X \times G$ with the same notations, then $iso(\mathcal G) = \{(x,x,g), x \in X, g \in G \}$, so any $h \in iso(\mathcal G)$ not in $G^0$, there is  $x \in X$ and  $g \in G$, g different from the unit element,   such that $h = (x,x,g)$, then $\alpha_h = \beta_g$ which can be taken in $Out(R)$ and obviouly $Z(A) = b(C(\mathcal G^0))$, the proposition follows from lemma \ref{explicit}.
\end{dm}

\subsection{Double crossed products}

Now, let's give a refinement of Proposition \ref{reconquete2} in the commutative case, that is $S = C(\mathcal G)$ (see also [Y] theorem 6.4 for  more general groupoids).

\subsubsection{\bf{Lemma}}
\label{coeur}
{\it   Let $\delta$ be a right action of  a commutative quantum groupoid $(C(\mathcal G), \Gamma_{\mathcal G},\kappa_{\mathcal G})$ on a von Neumann module $(A,b)$, then $\delta(A)( 1 \otimes C(\mathcal G) ) = \delta(1)(A \otimes C(\mathcal G))$.  }
\newline
\begin{dm}
Clearly, one has: $\delta(A)( 1 \otimes C(\mathcal G) \subset \delta(1)(A \otimes C(\mathcal G))$. On the other hand, using the identification of $ \delta(1)(A \otimes C(\mathcal G))$ with the set of functions $ \phi: \mathcal G \to A$ such that for any $g \in \mathcal G$, one has: $\phi(g) \in A_{r(g)}$, one easily sees that $\phi = \underset g \sum \delta(\alpha_{g^{-1}}(\phi(g)))(1 \otimes \delta_g)$, the lemma follows.
\end{dm}

\subsubsection{\bf{Theorem}}
\label{corazon}
{\it   Let $\delta$ be a right action of  a commutative quantum groupoid $(C(\mathcal G), \Gamma_{\mathcal G})$ on a von Neumann module $(A,b)$, then the double crossed product $ (A \underset \delta \rtimes C(\mathcal G))  \underset {\hat \delta} {\hat \rtimes} \mathcal R(\mathcal G) $ is isomorphic to $\delta(1)(A \otimes \mathcal W(\mathcal G))\delta(1)$, where $\mathcal W(\mathcal G)$ (the Weyl algebra of $C(\mathcal G$))  is the commutant in $\mathcal L( l^2(\mathcal G))$ of $r(l^\infty(\mathcal G^0))$ ($= \hat \alpha(N)'$) which is also the sub von Neumann algebra generated by $C(\mathcal G)$ and $\mathcal R(\mathcal G)$ ($=S\hat S$).  }
\vskip 1cm
\begin{dm}
Using Proposition \ref{reconquete2}, and the fact that $S = S'$, one deduces that: $ (A \underset \delta \rtimes C(\mathcal G))  \underset {\hat \delta} {\hat \rtimes} \mathcal R(\mathcal G) $ is isomorphic to $\delta(A)(1 \otimes S\hat S)\delta(1)$. Now thanks to lemma \ref{coeur} and the fact that $S\hat S$ is the vector space generated by $\{ s\hat s/ s \in S, \hat s \in \hat S \}$ (and also equal to $\beta(N)'$), the double crossed product is also isomorphic to $\delta(1)(A \otimes \mathcal S\hat S)\delta(1)$, the theorem follows.
\end{dm}

\subsection{Action of a groupoid on a fibered space over its base. }
\vskip 1cm

Let's suppose now that $A$ is commutative and finite dimensional, hence there is a finite set $X$ such that  $A= C(X)$, the von Neumann algebra of fonctions on $X$, the existence of $b$ leads to a partition $X= \underset{g \in \mathcal G^0} \sqcup X_u$, and for each $u \in \mathcal G^0$, one has: $A_u = C(X_u)$.

A left (resp. right) action of the groupoid $\mathcal G$ on $(A,b)$ is given by  a covariant (resp.contravariant) functor between the small category $\mathcal G$ and the category of sets with usual applications $\{X^u / \  u \in \mathcal G^0 \}$. So for any $g \in \mathcal G$, there exists an application  $ g \triangleright {\bf .}: X^{s(g)} \mapsto X^{r(g)}$ (resp $ {\bf . }\triangleleft g : X^{r(g)} \mapsto X^{s(g)}$)   such that   for any $g,g' \in \mathcal G$ which are composable, one has for any $x$ in $X^{s(g')}$ (resp. $X^{r(g)}$: $(g  \triangleright (g'  \triangleright x) = gg' \triangleright x$ (resp.$( x \triangleleft g) \triangleleft g' =  x\triangleleft g'g$). The bijection between the two notions is given by the following fomulae: 
\vskip 0.5cm
\centerline{for any $a$ in $A$ and $g$ in $\mathcal G$:  $\gamma(a)_g = g^{-1} \triangleright a_{r(g)}$ (resp. $\delta(a)_g = a_{s(g)} \triangleleft g$).}

\vskip 1cm

The crossed product  $C(\mathcal G) \underset \gamma \ltimes A $ (resp. $A \underset \delta \rtimes C(\mathcal G)$) can also be interpreted as the image of  a certain $*$-algebra representation. 

Let's denote by  $X_b \times _r\mathcal G = \{ (x,g) \in X \times \mathcal G / b(x) = r(g) \} $, the fiber product of $X$ and $\mathcal G$,  and by $L^1(X_b \times _r\mathcal G )$ the vector space of fonctions on this set. One can give to $L^1(X_b \times _r\mathcal G )$ a $*$-algebra structure denoted by $(L^1(X_b \times _r\mathcal G ), \underset \gamma \star , ^{ \underset \gamma \#})$ (resp.$(L^1(X_b \times _r\mathcal G ), \underset \delta \star , ^{ \underset \delta \#})$). 

For any fonctions $F,F'$ in  $L^1(X_b \times _r\mathcal G )$ and any $(x,g)$ in $X_b \times _r\mathcal G$, one has:

$$ F \underset \gamma \star F'(x,g)  =
 \underset{r(h) = r(g)}\sum F(x,h)F'(h^{-1}\triangleright x, h^{-1}g)$$
  
  $$ F^{ \underset \gamma \#} (x,g)  =
\overline{F(g^{-1}\triangleright  x, g^{-1}})$$

 $$\big( resp.  F \underset \delta \star F'(x,g)  =
 \underset{r(h) = r(g)}\sum F(x,h)F'(x \triangleleft h,  h^{-1}g)  $$
 
 $$ F^{ \underset \delta \#} (x,g)  =
\overline{F(x \triangleleft g, g^{-1}})  \  \  \big )$$
 One must keep in mind that these fonctions have the good support. One can define  a left (resp.right) regular representation of $L^1(X_b \times _r\mathcal G )$ in $l^2( X_b \times _s\mathcal G )$ (resp.  $l^2(X_b \times _r\mathcal G)$) denoted $L^\gamma$ (resp. $R^\delta$); for any $\xi$ in $l^2(X_b \times _s\mathcal G )$ (resp $l^2(X_b \times _r\mathcal G )$) any $F$ in $L^1(X_b \times _r\mathcal G )$ and any $(x,g)$ in $X_b \times _s\mathcal G$ (resp.$X_b \times _r\mathcal G$):
 
 $$ L^ \gamma(F)\xi(x,g) = \underset{r(h) = r(g)} \sum F(g \triangleright x, h) \xi(x, h^{-1}g)  $$
 
  $$\big( resp.  R^ \delta(F)\xi(x,g) = \underset{r(h) = s(g)} \sum F(x \triangleleft g, h) \xi(x, gh)  \big)$$

With these definitions one can also formulate an alternative definition of the crossed products:
$L^\gamma(L^1(X_b \times _r\mathcal G )) = C(\mathcal G )\underset \gamma \ltimes A $ and $R^\delta(L^1(X_b \times _r\mathcal G )) = A \underset \delta \rtimes C(\mathcal G)$

\section{Quantum groupoids coming from match  pairs of groupoids}

\subsection{The match pair of groupoids situation}
\label{match}
Now let's explain an  extension of the commutative example. Let $\mathcal G$ be any groupoid and $\mathcal H,\mathcal K$ be two subgroupoids of $\mathcal G$ such that $\mathcal G = \mathcal H\mathcal K =\{hk / h \in \mathcal H, k \in \mathcal K^{s(h)}\}$ and such that $\mathcal H \cap \mathcal K  \subset \mathcal G^0$, such a pair $\mathcal H,\mathcal K$ is called a match pair of groupoids (see [AA] for an abstract point of view). One easily verifies that this implies that $\mathcal G^0= \mathcal H \cap \mathcal K$ and that for any $g$ in $\mathcal G$ the decomposition $g = hk, h \in \mathcal H , k \in \mathcal K$ is unique. Hence one can  define two applications,  $p_1: \mathcal G \to \mathcal H$ and $p_2:  \mathcal G \to \mathcal K$ by the relation $g = p_1(g)p_2(g)$ for any $g \in \mathcal G$. Clearly one has $s \circ p_2 = s$ and $r \circ p_1 = r$, but a new application appears, the {\bf{middle}} one : 

\subsubsection{\bf Notations }
\label{middle}
{\it 
1) One has  $s \circ p_1 = r \circ p_2$, this application  will be denoted $\bf m$. 

2) For any  $f \in C(\mathcal G^0)$, we define $\alpha, \beta, \hat \beta$ by:
$\alpha(f) = f \circ m$ (the middle representation), $\beta(f) = f \circ r$ (the range representation) and $\hat \beta(f) = f \circ s$ (the source representation).

3)  With the exception of  the four representations of the base $N= C(\mathcal G^0)$, we shall use the same notations than in the   commutative case.  }

\subsubsection{\bf Lemma}
\label{trombone}
{\it For any $h$ in $\mathcal H$, $Card (\mathcal K^{s(h)}) = Card(\mathcal K^{r(h)})$ }
\newline
\begin{dm} 
Let's fix  $h$ in $\mathcal H$; let $k$ be any element of $\mathcal K^{s(h)}$. As $\mathcal G = \mathcal H\mathcal K$, then also $\mathcal G= \mathcal K\mathcal H$, so there exists a single pair $(k',h')$ in $\mathcal K\mathcal H$ such that $hk = k'h'$. Let's prove that the application $k \mapsto k'$ defines an injection from $\mathcal K^{s(h)}$ into $\mathcal K^{r(h)}$; if $k_1$ is any element of $\mathcal K^{r(h)}$ such that t$k'_1 = k'$, then there exists $h'_1$ for which one has  $hk_1 = k'h'_1$, one deduces that $h^{-1}k' = k_1h_1^{-1} = kh'^{-1}$ from which one deduces that $k= k_1$ and $h' = h'_1$. So the application is injective and $Card (\mathcal K^{s(h)}) \leq Card(\mathcal K^{r(h)})$, applying this to $h^{-1}$ one also has the inverse inequality. The lemma follows.
\end{dm}

\subsubsection{\bf Lemma}
\label{isotr}
{\it For any $u$ in $\mathcal G^0$, one has the following equalities: $Card(m^{-1}(u)) = Card(s^{-1}(u))  = Card(r^{-1}(u)) ( = Card (\mathcal G^u) = Card (\mathcal G_u))$. One has : $ tr \circ \alpha = tr \circ \beta = tr \circ \hat \beta. $ }
\newline
\begin{dm}
 The equality $Card(s^{-1}(u))  = Card(r^{-1}(u))$ is well known and is due to the bijection $g \mapsto g^{-1}$  which gives $ tr \circ \beta = tr \circ \hat \beta$. For every $f$ in $C(\mathcal G^0)$ one has:  $(tr \circ m)(f) = \sum_{u \in \mathcal G^0} Card(m^{-1}(u))f(u)$,  so the only thing to prove is that for any $u \in \mathcal G^0$, one has: $Card (\mathcal G^u)= Card(m^{-1}(u))$. 
 
But the application $(h,k) \mapsto hk$ is a bijection between $\mathcal H_u \times \mathcal K^u$ and $m^{-1}(u)$, so $$Card(m^{-1}(u)) = Card (\mathcal H_u)Card (\mathcal K^u ). $$
In an other hand, any element $g$ in $\mathcal G^u$ has a unique decomposition $g= hk$ where $h \in \mathcal H^u$ and $k \in \mathcal K^{s(h)}$, one easily gets that the image of $\mathcal G^u$ by the bijection  $g \mapsto (h,k)$ is equal to  the disjoint union: $\underset{h \in \mathcal H^u}\bigsqcup \{h\} \times \mathcal K^{s(h)}$, so using lemma \ref{trombone} and the last equality, one has:
\begin{align*}
 Card (\mathcal G^u) 
 &= \underset {h \in \mathcal H^u} \sum Card( \mathcal K^{s(h)}) =  \underset {h \in \mathcal H^u} \sum Card( \mathcal K^{r(h)}) = Card(\mathcal H^u)Card( \mathcal K^u) \\
 &= Card(\mathcal H_u)Card( \mathcal K^u) = Card (m^{-1}(u))
 \end{align*}
 \end{dm}

\subsubsection{\bf Lemma}
\label{metr}
{\it For any $x,y$ in $\mathcal G$ such that $m(x) = r(y)$, one has:

1) the elements $p_2(x)^{-1}$ and $y$  are composable for the multiplication of $\mathcal G$

2) the same is true for  $x$ and $p_1(p_2(x)^{-1}y)$,

3) $m(xp_1(p_2(x)^{-1}y)) = m(y)$. }
\newline
\begin{dm} For any $x,y$ in $\mathcal G$ such that $m(x) = r(y)$, then 
$ s(p_2(x)^{-1}) = r(p_2(x)) = m(x) = r(y)$, so $p_2(x)^{-1}$ and $y$ are composable. But one has: $r(p_1(p_2(x)^{-1}y) = r(p_2(x)^{-1}y) = r(p_2(x)^{-1}) = s(p_2(x)) = m(x) = s(x)$, so $x$ and $p_1(p_2(x)^{-1}y)$ are composable too. As for any  $(h,k)$ in $\mathcal H\times \mathcal K$ and $g$ in $\mathcal G$, one has: $m(hgk) = m(g)$, one deduces that :
$m(xp_1(p_2(x)^{-1}y)) = m(p_2(x)p_1(p_2(x)^{-1}y))$; let $( h_1, k_1)$ be in $\mathcal H \times \mathcal K$ and such that: $p_2(x)^{-1}y = h_1k_1$, then: $m(xp_1(p_2(x)^{-1}y)) = m(p_2(x)h_1) = m(yk_1^{-1}) = m(y)$.
\end{dm}

So the following definition is relevant:

\subsubsection{\bf Definition }
\label{belex}
{\it We shall denote $I_{\mathcal H,\mathcal K}$ the linear endomorphism of $l^2(\mathcal G)$ defined for any $f$ in $l^2(\mathcal G)$ and $x,y$ in $\mathcal G$ by:
$$  I_{\mathcal H,\mathcal K}(f)(x,y) = \left \{ \begin{array}{rl}
& f(xp_1(p_2(x)^{-1}y), p_2(x)^{-1}y) \hskip 0.5cm  \mathrm{if} \   m(x) = r(y) \\
&0 \hskip 0.5cm  \ \mathrm{otherwise}
\end{array}
\right.
$$
In particular: $I_{\mathcal G,\mathcal G^0} = I_\mathcal G$ and $I_{\mathcal G^0,\mathcal G}$ is the mpi studied in [Val1] 4.1. }

\subsubsection{\bf Proposition }
\label{cbo}
{\it $I_{\mathcal H,\mathcal K}$ is a mpi over the base $(C(\mathcal G^0), \alpha, \beta, \hat \beta)$. }
\newline
\begin{dm}
An easy computation gives the following formula for $\I^*$,  for any $f$ in $l^2(\mathcal G)$ and $x,y$ in $\mathcal G$ one has:
$$  I_{\mathcal H,\mathcal K}^*(f)(x,y) = \left \{ \begin{array}{rl}
& f(xp_1(y)^{-1}, p_2(xp_1(y)^{-1})y) \hskip 0.5cm  \mathrm{if} \   s(x) = m(y) \\
&0 \hskip 0.5cm  \ \mathrm{otherwise}
\end{array}
\right.
$$
So one has:
$$  I_{\mathcal H,\mathcal K}^* I_{\mathcal H,\mathcal K}(f)(x,y) = \left \{ \begin{array}{rl}
& f(x,y) \hskip 0.5cm  \mathrm{if} \   s(x) = m(y) \\
&0 \hskip 0.5cm  \ \mathrm{otherwise}
\end{array}
\right.
$$
and:
$$  I_{\mathcal H,\mathcal K}I_{\mathcal H,\mathcal K}^* (f)(x,y) = \left \{ \begin{array}{rl}
& f(x,y) \hskip 0.5cm  \mathrm{if} \   m(x) = r(y) \\
&0 \hskip 0.5cm  \ \mathrm{otherwise}
\end{array}
\right.
$$
This means  $I_{\mathcal H,\mathcal K}$ is a partial isometry the initial (resp.final) support of which is $e_{s,m}$ (resp. $e_{m,r}$). Let $f,f' $ be any element in  $C(\mathcal G^0)$ , $\xi$ be any element  in $l^2(\mathcal G\times \mathcal G)$,  $x,y$ be any element  in $\mathcal G$. 

First suppose  that $s(x) = m(y)$ then:
\begin{align*}
\I(\beta(f) \otimes \hat \beta(f')) &\xi(x,y) 
 = (\beta(f) \otimes \hat \beta(f'))\xi(xp_1(p_2(x)^{-1}y), p_2(x)^{-1}y) \\
&= f(r(xp_1(p_2(x)^{-1}y)f'(s( p_2(x)^{-1}y)\xi(xp_1(p_2(x)^{-1}y), p_2(x)^{-1}y)\\
&= f(r(x))f'(s( y))\xi(xp_1(p_2(x)^{-1}y), p_2(x)^{-1}y) \\
&= f(r(x))f'(s( y))\I \xi(x, y)\\
&= (\beta(f) \otimes \hat \beta(f'))\I\xi(x,y) 
\end{align*}
and, using lemma \ref{metr} 3):
\begin{align*}
\I(\alpha(f) \otimes  \beta(f')) & \xi(x,y) 
=(\alpha(f) \otimes \beta(f'))\xi(xp_1(p_2(x)^{-1}y), p_2(x)^{-1}y) \\
&= f(m((xp_1(p_2(x)^{-1}y)f'(r( p_2(x)^{-1}y)\xi(xp_1(p_2(x)^{-1}y), p_2(x)^{-1}y)\\
&= f(m(y))f'(r( p_2(x)^{-1})\xi(xp_1(p_2(x)^{-1}y), p_2(x)^{-1}y) \\
&= f(m(y))f'(s(x))\I \xi(x, y)\\
&= (\hat \beta(f') \otimes \alpha(f))\I\xi(x,y) 
\end{align*}
If  $s(x) \not= m(y)$, one has: $\I(\beta(f) \otimes \hat \beta(f'))\xi(x,y)  = 0 = (\beta(f) \otimes \hat \beta(f'))\I\xi(x,y)  $, and also: $\I(\alpha(f) \otimes  \beta(f'))\xi(x,y) = 0 = (\hat \beta(f') \otimes \alpha(f))\I\xi(x,y) $. Hence: $$\I(\beta(f) \otimes \hat \beta(f')) = (\beta(f) \otimes \hat \beta(f'))\I $$
$$ \I(\alpha(f) \otimes  \beta(f'))  =  (\hat \beta(f') \otimes \alpha(f))\I$$

Now let's prove the pentagonal relation for $\I$. Let's fix some notation:  for any $x,y$ in $\mathcal G$ such that $m(x) = r(y)$ then one can define: $V = p_2(x)^{-1}y$ and  $X=  xp_1(V)$,  if  moreover $z$ is any element of $\mathcal G$ such that $m(y) = r(z)$, then $m(p_2(x)^{-1}y) = r(p_2(X)^{-1}z)$ and, by two routine calculations  the following relations are true:

\begin{multline*}
(\I)_{12}(\I)_{13}(\I)_{23}\xi(x,y,z)= \\
=\xi(xp_1(V)p_1(p_2(X)^{-1}z),p_2(x)^{-1}yp_1(p_2(V)^{-1}p_2(X)^{-1}z),p_2(V)^{-1}p_2(X)^{-1}z)
\end{multline*}
\begin{multline*}
(\I)_{23}(\I)_{12}\xi(x,y,z)= \\
=\xi(xp_1(p_2(x)^{-1}yp_1(p_2(y)^{-1}z)),p_2(x)^{-1}yp_1(p_2(y)^{-1}z),p_2(y)^{-1}z)
\end{multline*}

Let $(h,k)$ be in $\mathcal H \times \mathcal K$ such that: $p_2(x)^{-1} p_1(y) = hk$, then $V= hkp_2(y)$ and $X= p_1(x)p_2(x)h = p_1(x)p_1(y)k^{-1}$ hence one has:

\begin{align}
&p_2(V)^{-1}p_2(X)^{-1}= p_2(y)^{-1}  \\
&p_2(X)^{-1} = k  \\
&p_1(V) = h
\end{align}

So, using (10) and the notation:  $k' = kp_2(y)$, one has:
\begin{align*}
&p_1(p_2(x)^{-1}yp_1(p_2(y)^{-1}z))=\\
&= p_1(hkp_2(y)p_1(p_2(y)^{-1}z))=  hp_1(kp_2(y)p_1(p_2(y)^{-1}z)) =  p_1(V)p_1(k'p_1(k'^{-1}kz))
\end{align*}

Now let's define $(h' , k'')$ in $\mathcal H\times \mathcal K$ such that: $k'^{-1}kz = h'k''$, hence using (9), one has:
\begin{align*}
p_1(p_2(x)^{-1}y&p_1(p_2(y)^{-1}z))= \\
&=  p_1(V)p_1(k'h') =  p_1(V)p_1( kzk'')  = p_1(V)p_1( kz) = p_1(V)p_1( p_2(X)^{-1}z)
\end{align*}

This last equality and (8) gives that for any triple $(x,y,z)$ in $\mathcal G^3$ such that $m(x) = r(y)$ and $m(y) = r(z)$:
$$ (\I)_{12}(\I)_{13}(\I)_{23}\xi(x,y,z)= (\I)_{23}(\I)_{12}\xi(x,y,z)$$
but for all the other triples $(x,y,z)$ in $\mathcal G^3$ the two sides of this equality are $0$, hence, $\I$ is a mpi.
\end{dm}

\subsection{Crossed products and match pairs of groupoids}
The situation of a match pair of groupoids $\mathcal G =\mathcal H\mathcal K$, leads to a natural right action of the groupoid $\mathcal H$ on the fibered space  $\mathcal K$ and a left action of the groupoid $\mathcal K$ on the fibered space  $\mathcal H$. Using the inverse map, one has $\mathcal G = \mathcal H\mathcal K = \mathcal K\mathcal H$.  Hence, for any $k \in \mathcal K$ and $h\in \mathcal H^{s(k)}$,  there exist a unique $h' \in H$ and a unique $k' \in K^{s(h')}$ such that $kh = h'k'$.

\subsubsection{\bf{Lemma and definition}}
\label{bienfait}
{\it  Let $\mathcal G =\mathcal H\mathcal K$ be a match pair of groupoids, and for any   $k \in \mathcal K$ and $h\in \mathcal H^{s(k)}$, let's denote by $k \triangleright h $ (resp. $k \triangleleft h$) the unique element in $H$ (resp. $K^{s(k \triangleright h )} )$ such that:
 $$kh= (k \triangleright h)(k \triangleleft h),$$
then $ \triangleright$ (resp.$\ \triangleleft$) is a left action of the groupoid $\mathcal K$ on the fibered space  $\mathcal H$ (resp. right action of the groupoid $\mathcal H$ on the fibered space  $\mathcal K$) }
\newline
\begin{dm} Left to the reader
\end{dm}

Let us denote by $\mathcal G/\mathcal K$ (resp $\mathcal H\backslash \mathcal G$) the set of right (resp.left) classes in $\mathcal G$ modulo $\mathcal K$ (resp.$\mathcal H$), that is $\{ g\mathcal K^{s(g)} /  \  g \in \mathcal G\}$ (resp. $\{\mathcal H_{r(g)}g /  \  g \in \mathcal G\}$). In that case, the application : $h \to h\mathcal K^{s(h)}$ (resp. $k \to \mathcal H_{r(k)}k$) is  a natural bijection between  $\mathcal H$ and $\mathcal G/\mathcal K$ (resp. $\mathcal K$ and $\mathcal H\backslash \mathcal G$). Using these applications, $\mathcal G/\mathcal K$ and $\mathcal H\backslash \mathcal G$ are fibered by $\mathcal G^0$: for any $u$ in $\mathcal G^0$, one can define $(\mathcal G/\mathcal K)^u = \{ g\mathcal K^{s(g)} /  \ r(g) = u \}$ and $(\mathcal H\backslash \mathcal G)^u = \{ \mathcal H_{r(g)}g/ \ s(g) = u \}$. Also $\mathcal K$ (resp. $\mathcal H$) has a left action on  $\mathcal G/\mathcal K$  (resp.  $\mathcal H\backslash \mathcal G$) by multiplication: for any $h$ in $\mathcal H$, $k$ in $\mathcal K$, $g$ in $\mathcal G_{r(h)}$ and $g'$ in $\mathcal G^{s(k)}$, one can define $\delta_h(\mathcal H_{r(g)}g) = \mathcal H_{r(g)}gh$ and $\gamma_k(g'\mathcal K^{s(g')}) = kg'\mathcal K^{s(g')}$. Using the natural bijections below, and slightly abusing notations, one easily sees that the right action  of $ \mathcal H$ on $\mathcal K$ (resp. left action of $\mathcal K$ on $\mathcal H$) is exactly the one coming from lemma \ref{bienfait}

Now, in the following, let's remember that $\chi_p$ (resp $\chi_Z$) denotes the characteristic fonction of the singleton $\{p\}$ (resp. set $Z$).

\subsubsection{\bf{Proposition and notations}}  
\label{dur}
 {\it   The mpi $\I$ is regular and the $C^*$-algebra $S$ (resp.$\hat S$) associated to $\I$ is isomorphic to the crossed product $C(\mathcal K )\underset \k \ltimes C(\mathcal H)$ (resp. $ C(\mathcal K) \underset \h \rtimes C(\mathcal H)$), where $\gamma^{\mathcal K}$ is the left   action of $C(\mathcal  K)$ on $C(\mathcal H)$ associated with $\triangleright$ (resp.$ \delta^{\mathcal H}$ is the right  action of $C(\mathcal  H)$ on $C(\mathcal K)$ associated with $\triangleleft$). Hence $ C(\mathcal K) \underset {\delta^{\mathcal H} }\rtimes C(\mathcal H)$ and $C(\mathcal K) \underset {\gamma^{\mathcal K}} \ltimes C(\mathcal H)$ have  weak Hopf $C^\star$-algebras structures in duality, we shall note them $(C(\mathcal K )\underset {\gamma^{\mathcal K}} \ltimes C(\mathcal H), \Gamma^\gamma, \kappa^\gamma, \epsilon^\gamma)$ and $(C(\mathcal K) \underset {\delta^{\mathcal H} }\rtimes C(\mathcal H), \Gamma^\delta, \kappa^\delta, \epsilon^\delta)$.  }
\vskip 1cm
\begin{dm}
For any $g,g',p,q$ in $\mathcal G$ and $\xi$ in $l^2(\mathcal G)$, one has:
\begin{align*}
(i\otimes \omega_{\chi_p,\chi_q})(\I)(\xi)(g)
&= ((i\otimes \omega_{\chi_p,\chi_q})(\I)(\xi), \chi_g) = (\omega_{\xi,\chi_g} \otimes \omega_{\chi_p,\chi_q})(\I) \\
&= (\I(\xi \otimes \chi_p), \chi_g \otimes \chi_q) = \I(\xi \otimes \chi_p)(g,q) \
\end{align*}
Hence,  $(i\otimes \omega_{\chi_p,\chi_q})(\I)\xi(g) = 0$ if $m(g) \not = r(q)$;  otherwise: 
\begin{align*}
(i\otimes \omega_{\chi_p,\chi_q})(\I)(\xi)(g)
&= \I(\xi \otimes \chi_p)(g,q) = \chi_p(p_2(g)^{-1}q)\xi(gp_1(p_2(g)^{-1}q))
\end{align*} 

This also implies that $(i\otimes \omega_{\chi_p,\chi_q})(\I)(\xi)(g) \not = 0$ only if there exists $k \in \mathcal K$ such that $q = kp$ and $p_2(g) = k$; one can see that these two conditions imply that $m(g) = r(q)$. 

To resume,  for any $k$ in $\mathcal K$, $p$ in $\mathcal G^{s(k)}$  and $ g $ in $\mathcal G$, one has:
$(i\otimes \omega_{\chi_{p},\chi_{kp}})(\I)(\xi)(g) = 1_{p_2^{-1}(k)}(g)\xi(gp_1(p))$, and if $q$ is  not in $\mathcal Kp$, one has: $(i\otimes \omega_{\chi_{p},\chi_{kp}})(\I)= 0$. So $\hat S$ is generated by the operators $ : \xi \mapsto (g \mapsto  \chi_{\mathcal Hk}(g)\xi(gh))$, for any $(k,h)$ in $\mathcal K_s  \times  _r\mathcal H$, up to the natural identification of  $ \mathcal G$ with $\mathcal K_s  \times  _r\mathcal H$, this is $ R^{\delta}(\chi_{(k,h)})$ so $\hat S$ is isomorphic to the crossed product $ C(\mathcal K) \underset{\delta^{\mathcal H}} \rtimes C(\mathcal H)$, hence $\I$  is regular. In a very similar way,   $S$ is generated by the operators $(\omega_{\chi_{kh},\chi_k} \otimes i)(I_{\mathcal H, \mathcal K})$ for $(h,k)$ in $\mathcal K_s  \times  _r\mathcal H$ which appear, up to the identification of $\mathcal H_r \otimes _r\mathcal K$ with $\mathcal G$, to be equal to $L^\gamma(\chi_{(k \triangleright h, k)})$ (observe that $r(k) = r(k \triangleright h)$); hence $S$ is isomorphic to $C(\mathcal K) \underset {\gamma^{\mathcal K}} \ltimes C(\mathcal H)$ .
\end{dm}
\vskip 0.5cm

Let's  compare these structures to  ones defined in by N. Andruskiewitsch  and S.Natale [AN]. Let's use the notations of [AN] theorem 3.1, and let's identify $\mathcal K_s  \times  _r\mathcal H$ and $\mathcal T$, the double groupoid associated with the match pair $\mathcal K\mathcal H$ by [AN] proposition 2.9, using the bijection:
 $$(k,h) \mapsto  \eenmatrix{k}{h}{ }{ }$$

As vector spaces $ C(\mathcal K) \underset{\delta^{\mathcal H}} \rtimes C(\mathcal H) = \mathbb C \mathcal T$. As one can define on $\mathcal T$ an horizontal and a vertical product (see lemma 1.5 of [AN]),  the above identification gives rise to compositions laws on $\mathcal K_s  \times  _r\mathcal H$.

\subsubsection{\bf Definition and notations }
\label{couscous}
{\it  We shall denote by $\overset \rightarrow \square$  the horizontal  product defined for every $(k,h), (k',h')$ in  $\mathcal K_s \times _r \mathcal H$ such that $h' = k \triangleright h$ by :
$$(k,h) \overset \rightarrow \square (k', h') = (kk',h').$$
We shall denote by $\square \downarrow$ the vertical product defined for every $(k,h), (k',h')$ in  $\mathcal K_s \times _r \mathcal H$ such that $k' = k \triangleleft h$ by :
$$(k,h) \square \downarrow(k', h') = (k,hh').$$}

Now we shall give a complete description of the  $C^*$-quantum groupoid structure given by proposition \ref{dur} to $ C(\mathcal K) \underset{\delta^{\mathcal H}} \rtimes C( \mathcal H)$, which proves that it is isomorphic to
 $\mathbb C\mathcal T$.

\subsubsection{\bf Lemma }
\label{couscous1}
{\it One has: $ \Gamma^\delta \delta^{\mathcal H}= (\delta^{\mathcal H}\otimes \delta^{\mathcal H})\Gamma_{\mathcal K}$, so for any $k$ in $\mathcal K$ and any $g,g'$ in $\mathcal G$, one gets: 
$$ \Gamma^\delta (\delta^{\mathcal H}(\chi_k))\xi(g,g') = \chi_k(p_2(g)p_2(g'))\xi(g,g')$$}
\newline
\begin{dm}
Left to the courageous reader.
\end{dm}

\subsubsection{\bf Theorem }
\label{pantalon}
{\it 1) For every $(k,h), (k',h')$ in  $\mathcal K_s \times _r \mathcal H$, one has: 

$  \chi_{(k,h)}) \underset \delta \star \chi_{(k',h')}  = \left \{ \begin{array}{rl}
&\chi_{(k,h) \square \downarrow (k',h')})   \hskip 0.5cm  \mathrm{if}  \   k' = k \triangleleft  h\\
&0 \hskip 0.5cm  \ \mathrm{otherwise}
\end{array}
\right.$ \ , \  and $ \chi_{(k,h)})^{ \underset \delta \#} = \chi_{(k \triangleleft  h, h^{-1})}$

The identification of $\mathcal K_s \times _r \mathcal H$ with $\mathcal T$ given by the application: $(k,h) \mapsto  \eenmatrix{k}{h}{ }{ }$, leads to a  $C^*$-isomorphism between   and $C(\mathcal K) \underset {\delta^{\mathcal H} }\rtimes C( \mathcal H) $ and $\mathbb C \mathcal T$.

2) For every $(k,h)$ in  $\mathcal K_s \times _r \mathcal H$, one has: 

$$ \Gamma^\delta( R^{\delta}(\chi_{(k,h)})) =  \underset{(k_1,h_1)  \overset \rightarrow \square (k_2, h_2) = (k,h)}\sum R^{\delta}(\chi_{(k_1,h_1)} ) \otimes R^{\delta}(\chi_{(k_2,h_2)}).$$ 

$$  \kappa^\delta(R^{\delta}(\chi_{(k,h)})) = R^\delta(\chi_{((k \triangleleft h)^{-1}, (k \triangleright h)^{-1})})$$

$$ \epsilon(R^\delta(\chi_{(k,h)})) =   \left \{ \begin{array}{l}
1  \  \   \mathrm{if}  \  k  = r(h) )\\
0  \  \   \mathrm{otherwise} 
\end{array}
\right.$$

$R^{\delta}$ is an isomorphism of $C^*$-quantum groupoids between $\mathbb C \mathcal T$ associated with the match pair  $\mathcal K \mathcal H$ by proposition 3.4 of [AN] and $(C(\mathcal K) \underset {\delta^{\mathcal H} }\rtimes C(\mathcal H), \Gamma^\delta, \kappa^\delta, \epsilon^\delta)$. }
\vskip 1.5cm
\begin{dm}
An easy computation gives that, for all $(k,h), (k',h')$ in  $\mathcal K_s \times _r \mathcal H$, one has: 

$  \chi_{(k,h)}) \underset \delta \star \chi_{(k',h')}  = \left \{ \begin{array}{rl}
&\chi_{(k,h) \square \downarrow (k',h')})   \hskip 0.5cm  \mathrm{if}  \   k' = k \triangleleft  h\\
&0 \hskip 0.5cm  \ \mathrm{otherwise}
\end{array}
\right.$ \ , \  and $ \chi_{(k,h)})^{ \underset \delta \#} = \chi_{(k \triangleleft  h, h^{-1})}$, but $(k \triangleleft  h, h^{-1})$ appears to be the inverse of $(h,k)$ for the vertical product . 
\vskip 0.2cm
Hence, if one denotes by char($\eenmatrix{k}{h}{ }{ })$ the characteristic fonction of $\eenmatrix{k}{h}{ }{ }$ in $\mathbb C \mathcal T$, the  application: 

$$ char(\eenmatrix{k}{h}{ }{ } )\mapsto R^\delta(\chi_{(k,h)} ) $$

gives an explicit   $*$-algebras isomorphism between $ \mathbb C \mathcal T$ and  $ C(\mathcal K) \underset{\delta^{\mathcal H}} \rtimes C( \mathcal H)  $.

Let's denote by $\Theta$ the image by this isomorphism  of the coproduct $\Delta$ given by theorem 3.1 of [AN],  so for any $(k, h)$ in $\mathcal K_s \times _r \mathcal H$:  

 $ \Theta(R^{\delta}(\chi_{(k,h)}) ) =  \underset{(k_1,h_1)  \overset \rightarrow \square (k_2, h_2) = (k,h)}\sum R^{\delta}(\chi_{(k_1,h_1)} ) \otimes R^{\delta}(\chi_{(k_2,h_2)})$

Let's first prove the first formulae of the theorem when $h = s(k)$, ( i.e  $ \Gamma^\delta(\delta^{\mathcal H}(\chi_k) ) =  \Theta(\delta^{\mathcal H}(\chi_k) )$. One can easily observe that  for any $(k_1,h_1), (k_2,h_2)$ in $\mathcal K_s \times _r \mathcal H$, one has $(h_1,k_1)  \overset \rightarrow \square (k_2, h_2) = (k,s(k))$ if and only if $k_1k_2 = k$ and $h_1= s(k_1), h_2 = s(k_2)$, hence for any $g,g'$ in $\mathcal G$ and $\xi$ in $l^2(\mathcal G \times \mathcal G)$, by lemma \ref{couscous1} one has:
\begin{align*}
&
\underset{(k_1,h_1)  \overset \rightarrow \square (k_2, h_2) = (k,h)}\sum (R^{\delta}(\chi_{(k_1,h_1)} \otimes R^{\delta}(\chi_{(k_2,h_2)}))\xi(g,g') =  \\
&= \underset{k_1k_2= k}\sum (R^{\delta}(\chi_{(k_1,s(k_1))} \otimes R^{\delta}(\chi_{(k_2,s(k_2))}))\xi(g,g') \\
&=  \underset{k_1k_2= k}\sum  \  \  \underset{r(h') = s(g)} \sum  \  \  \underset{r(h'') = s(g')} \sum \chi_{(k_1,s(k_1))}(p_2(g),h') \chi_{(k_2,s(k_2)}(p_2(g'),h''))\xi(gh',g'h'') \\
&=  \underset{k_1k_2= k}\sum \chi_{k_1}(p_2(g)) \chi_{k_2}(p_2(g'))\xi(g,g') =  \chi_{k}(p_2(g)p_2(g'))\xi(g,g') \\
&= \Gamma^\delta (\delta^{\mathcal H}(\chi_k))\xi(g,g') = \Gamma^\delta (R^\delta(\chi_{(k,h)}))\xi(g,g').
\end{align*}

Let's denote that the fourth equality is due to the fact that:  $p_2(g) = k_1$ implies that $s(k_1) = s(g)$ and that $p_2(g') = k_2$ implies that $s(k_2) = s(g')$. So  one has: 

\begin{align}
\Gamma^\delta(\delta^{\mathcal H}(\chi_k) )  =\Theta(\delta^{\mathcal H}(\chi_k) )  \  \mathrm{for \  any }  \  k  \    Ê \mathrm{in}   \  \mathcal K
\end{align}

Now let's denote by $\rho $ the right regular representation of $\mathcal H$,  so for any $h,h'$ in $\mathcal H$ , $\xi$ in $l^2(\mathcal H)$:  $\rho(h)\xi(h') = \xi(h'h)$. Then in $ C(\mathcal K) \underset {\delta^{\mathcal H} }\rtimes C(\mathcal H)$, one has: $ 1 \otimes \rho(h) = \underset {s(k) = r(h)} \sum R^\delta(\chi_{(k,h)})$. As by definition, for any $h$ in $\mathcal H$, $\xi$ in $l^2(\mathcal G \times \mathcal G)$, $g,g' $ in $\mathcal G$, one has: $   \Gamma^\delta(1_{\mathcal K} \otimes \rho(h))\xi(g,g') =    I_{\mathcal H,\mathcal K}^*(1_{\mathcal G} \otimes ( 1_{\mathcal K} \otimes \rho(h))) I_{\mathcal H,\mathcal K}\xi(g,g') $,  an easy computation gives that:

$$   \Gamma^\delta(1_{\mathcal K} \otimes \rho(h))\xi(g,g') =\left \{ \begin{array}{rl}
&\xi(gp_1(p_2(g')h), g'h) \hskip 0.5cm  \mathrm{if}  \  s(g) = m(g') \    \mathrm{and} \   s(g') = r(h)\\
&0 \hskip 0.5cm  \ \mathrm{otherwise}
\end{array}
\right.
$$

In an other hand, for any  $k_1,k_2,K$ in $\mathcal K$,  $h_1,h_2,H$ in $\mathcal H$, one has $(h_1,k_1)  \overset \rightarrow \square (k_2, h_2) = (k,h)$ if and only if  $r(k_1) = r(k)$, $h_1 = k_1^{-1}k \triangleright h$, $k_2 = k_1^{-1}k$, $h_2 = h$. Hence, we can write that:

\begin{align*}
&\Theta(1_{\mathcal K} \otimes \rho(h))\xi(g,g') \\
&= \underset{s(k) = r(h)} \sum \   \underset {r(k_1) = r(k)}\sum  {R^\delta}(\chi_{(k_1,k_1^{-1}k \triangleright h)} ) \otimes R^\delta(\chi_{(k_1^{-1}k ,h)})\xi(g,g')  \\
&= \underset{s(k) = r(h)} \sum   \underset {r(k_1) = r(k)}\sum   \underset{r(h') = s(g)} \sum    \underset{r(h'') = s(g')} \sum  \chi_{k_1}(p_2(g))\chi_{k_1^{-1}k \triangleright h} (h') \chi_{k_1^{-1}k}(p_2(g'))\chi_{h} (h'')\xi(gh',g''')
\end{align*}

Due to the characteristic fonctions, all terms of the sum are zero except when:
$$  \left \{ \begin{array}{l}
k_1 = p_2(g)\\
k_2^{-1}k \triangleright h = h'\\
k_2^{-1}k  = p_2(g')\\
h'' = h
\end{array}
\right.
$$

In order to have non zero terms, one needs that $k = p_2(g)p_2(g')$, hence $s(g) = m(g')$, and also $s(g') = s(k)$,  which implies $s(g') = r(h)$. In these conditions there is a single term in the sum, it is  obtained for: $k = p_2(g)p_2(g'), k_1 = p_2(g), h' = p_2(g)  \triangleright h = p_1(p_2(g)h)$, hence: $\Theta(1_{\mathcal K} \otimes \rho(h))\xi(g,g') = \xi(gp_1(p_2(g')h), g'h)$. One deduces that:
\begin{align}
  \Gamma^\delta(1_{\mathcal K} \otimes \rho(h)) = \Theta(1_{\mathcal K} \otimes \rho(h)) 
\end{align}

As both $ \Gamma^\delta$ and $\Theta$ are multiplicative, using (4) and (5), one deduces that for any $(k, h)$ in $\mathcal K_s \times _r \mathcal H$:  
 $$ \Theta(R^{\delta}(\chi_{(k,h)}) ) =  \underset{(k_1,h_1)  \overset \rightarrow \square (k_2, h_2) = (k,h)}\sum R^{\delta}(\chi_{(k_1,h_1)} ) \otimes R^{\delta}(\chi_{(k_2,h_2)})$$
 
 Now for any $g$ in $\mathcal G$, any $(k, h)$ in $\mathcal K_s \times _r \mathcal H$, and any $\xi$ in $l^2(\mathcal G)$, one has:
 \begin{align*}
\kappa^\delta (R^\delta{\chi_{(k,h)}} )\xi(g) 
&= (i \otimes \omega_{\chi_h, \chi_{kh}})(I_{\mathcal H,\mathcal K}^*)\xi(g) 
= I_{\mathcal H,\mathcal K}^*(\xi \otimes \chi_h)(g, kh) \\
&= \xi (gp_1(kh)^{-1})\chi_h(p_2( gp_1(kh)^{-1})kh)
\end{align*}

One easily sees that $p_2( gp_1(kh)^{-1})kh = h$ if and only if $p_2(g) = (k \triangleleft h)^{-1}$

So: $\kappa^\delta (R^\delta{\chi_{(k,h)}} )\xi(g)  =  \xi (g( k \triangleright h)^{-1})\chi_{H(k \triangleleft h)}(g)= R^\delta(\chi_{((k \triangleleft h)^{-1}, (k \triangleright h)^{-1})})\xi(g)$

Hence:

$$  \kappa^\delta(R^{\delta}(\chi_{(k,h)})) = R^\delta(\chi_{((k \triangleleft h)^{-1}, (k \triangleright h)^{-1})})$$

For any $(k, h)$ in $\mathcal K_s \times _r \mathcal H$, one has:

$\epsilon(R^\delta(\chi_{(k,h)})) = \epsilon((i \otimes \omega_{\chi_h, \chi_{kh}})(I_{\mathcal H,\mathcal K}    ) =  \omega_{\chi_h, \chi_{kh}}(1) = \underset {g \in \mathcal G} \sum \chi_h(g) \overline \chi_{kh}(g)$

Hence:

$$ \epsilon(R^\delta(\chi_{(k,h)})) =   \left \{ \begin{array}{l}
1  \  \   \mathrm{if}  \  k  = r(h) )\\
0  \  \   \mathrm{otherwise} 
\end{array}
\right.$$
The theorem follows immediatly.
\end{dm}

\subsubsection{\bf{Remark }}
{\it Using the natural identification of $\mathcal G$ and $\mathcal K_s \times _r \mathcal H$, one can also express the $C^*$-quantum groupoid structure $(C(\mathcal K) \underset {\delta^{\mathcal H} }\rtimes C(\mathcal H), \Gamma^\delta, \kappa^\delta, \epsilon^\delta)$ by the following formulaes:

\begin{itemize}

\item 

$\Gamma^\delta( R^{\delta}(\chi_g)) =  \underset{(k_1,h_1)  \overset \rightarrow \square (k_2, h_2) = (p_2(g),p_1(g))}\sum R^{\delta}(\chi_{h_1k_1} ) \otimes R^{\delta}(\chi_{h_2k_2}).$

\vskip 0.2cm
\item

$\kappa^\delta(R^{\delta}(\chi_g)) = R^\delta(\chi_{g^{-1}})$
\vskip 0.3cm
\item
$ \epsilon(R^\delta(\chi_g)) =   \left \{ \begin{array}{l}
1  \  \   \mathrm{if}  \  g \in \mathcal H \\
0  \  \   \mathrm{otherwise} 
\end{array}
\right.$

\end{itemize}
}

Naturally, using the identification of  $\mathcal G$ and $\mathcal H_s \times _r \mathcal K$, one also has a characterization of the  $C^*$-quantum groupoid structure for $(C(\mathcal K) \underset {\gamma^{\mathcal K}} \ltimes C(\mathcal H), \Gamma^\gamma, \kappa^\gamma, \epsilon^\gamma)$ dual to $(C(\mathcal K) \underset {\delta^{\mathcal H} }\rtimes C(\mathcal H), \Gamma^\delta, \kappa^\delta, \epsilon^\delta)$ . Let's recall that $\mathcal T^t$ (the transpose of $\mathcal T$) is by definition equal to $\mathcal T$ as a set but the horizontal and vertical laws are exchanged and due to Proposition 3.11 [AN] $\mathbb C\mathcal T^t$ is the dual of $\mathbb C\mathcal T$. 
\vskip 0.5cm
The application $ (h,k) \mapsto  \eenmatrix{k}{ }{ h}{ }$ gives a  bijection between $\mathcal H_r \times _r \mathcal K$ and $\mathcal T^t$ (equal to $\mathcal T$ as a set).  

The same calculations than above give the following result:

\subsubsection{{\bf Proposition}}
\label{retour}
{\it There exists an isomorphism of $C^*$-quantum groupoids between $(C(\mathcal K ) \underset {\gamma^{\mathcal K}} \ltimes C(\mathcal H), \Gamma^\gamma, \kappa^\gamma, \epsilon^\gamma)$ and $\mathbb C \mathcal T^t$.}

One can give explicitly this isomorphism. If $char( \eenmatrix{h}{k}{ }{ })$ is the characteristic fonction of $ \eenmatrix{h}{k}{ }{ }$ in $\mathbb C \mathcal T^t$, then:
$$  L^\gamma(\chi_{(h,k)})  \mapsto char( \eenmatrix{k}{ }{h }{ }  )   $$.

 In [BH] is given a very deep study of inclusions of the form $R^H \subset R \rtimes K$ where $H$ and $K$ are subgroups of a group $G$ acting properly and outerly on the  hyperfinite type $II_1$ factor $R$, in such a way we can identify $G$ with a subgroup of $Out R$,  in particular there is here no ambiguity for $H \cap K \subset Out R$: let's call $\alpha$ the action of $K$ and $\beta$ the one of $H$  here these actions coincide on $H \cap K$ .    In [BH],   it's proved  that this inclusion is finite depth if and only if the group generated by $H$ and $K$ in $Out R$ is finite,  and, in that situation,  it is irreducible and depth two when $H,K$ is a match pair. Relaxing the match pair property, considering the  case $card(H \cap K) \not = 1$, let's prove now that the inclusion $R^H \subset R \rtimes K$ does not come from the match pair of groupoids procedure.

\subsubsection{\bf{Lemma(HB)} }
\label{dietmar}
{\it  In the preceeding conditions, the algebra $(R^H)' \cap R \rtimes K$ is isomorphic to 
the group algebra $\mathbb C[H \cap K]$. }
\newline
\begin{dm}
If $u_k$ and $v_h$ are canonical implementations of $\alpha$ and $\beta$ on $L^2(R)$, one can suppose $u_x = v_x$ for any $x$ in $H \cap K$, these $u_x$ generate a *-algebra isomorphic to $\mathbb C[H \cap K]$ and are clearly in $(R^H)' \cap R \rtimes K$, in an other hand, using the calculation in the proof of 4.1 in [HB], one has: $dim ((R^H)' \cap R \rtimes K )= card (H \cap K)$. The lemma follows.
\end{dm}

\subsubsection{\bf{Corollary}}
\label{photocopie}
{\it The algebra $(R^H)' \cap R \rtimes K$  is commutative if and only if $H \cap K$ is abelian. Hence, when $H \cap K$ is non abelian,  the inclusion $R^H \subset R \rtimes K$ does not come from the match pair of groupoids procedure.. }

\vskip 2cm
REFERENCES

[AA] M.AGUIAR \& N.ANDRUSKIEWITSCH, Representations of match pairs of groupoids and applications to weak hopf algebras {\it Contemp. Math.}  {\bf 376} (2005).

[AN] N.ANDRUSKIEWITSCH \& S.NATALE, \  \ Double categories and  quantum  groupoids {\it Publ. Mat. Urug. } {\bf 10} 11-51 (2005) .

[BS] S. BAAJ \& G. SKANDALIS, Unitaires multiplicatifs et dualit\'e pour 
les produits crois\'es de C*-alg\`ebres. {\it Ann. Sci. ENS} {\bf 26} 
(1993), 425-488.

[BBS] S. BAAJ  \& E. BLANCHARD \& G. SKANDALIS, Unitaires multiplicatifs 
en dimension finie et
leurs sous-objets. {\it Ann.Inst. Fourier} {\bf 49} (1999), 1305-1344.

[BH] D.BISCH \& U.HAAGERUP, Composition of subfactors: new examples of infinite depth subfactors, {\it Ann. Sci. ENS 4\`eme  s\' erie } {\bf 29} $n^o 3$ (1996), 329-383.

[BoSz] G. B\"OHM \& K.SZLACH\'ANYI, Weak C*-Hopf algebras: the  
coassociative symmetry of non integral dimensions, {\it  Quantum groups 
and quantum spaces. Banach Center Publications} {\bf 40} (1997), 9-19.

[BoSzNi] G. B\"OHM, K.SZLACH\'ANYI \& F.NILL, Weak Hopf Algebras I. 
Integral Theory and $C^*$-structure. {\it Journal of Algebra} {\bf 221} 
(1999),  385-438.

[E] M. ENOCK , Produit crois\'e d'une alg\`ebre de von Neumann par une alg\`ebre de Kac.
{\it Journal of Functional Analysis} {\bf 26} (1977), 16-47.

[E2] M. ENOCK , Inclusions of von Neumann algebras and quantum groupoids III .
{\it Journal of Functional Analysis} {\bf 223} (2005), 311-364.

[EV] M. ENOCK \& J.M. VALLIN, Inclusions of von Neumann
algebras and quantum groupoids.
{\it Journal of Functional Analysis} {\bf 172} (2000), 249-300.

[GHJ] F.M.GOODMAN, P. de la HARPE, V.F.R. JONES, Coxeter graphs and towers of algebras {\it M.S.R.I. publications} 14

[H] U. HAAGERUP, The standard form of von Neumann algebras.{\it Math Scand.}{\bf 37} (1975) 271-283

    [L] F.LESIEUR, thesis, http://tel.ccsd.cnrs.fr/documents/archives0/00/00/55/05 
 
[NV1] D. NIKSHYCH \& L. VAINERMAN, Algebraic versions of a 
finite-dimensional quantum groupoid. {\it Lecture Notes in Pure and 
Appl. Math.} {\bf 209} (2000), 189-221.

[NV2] D. NIKSHYCH \& L. VAINERMAN, A characterization of depth 2 
subfactors of $II_1$ factors. {\it JFA} {\bf 171} (2000), 278-307
(2000).

[NV3] D. NIKSHYCH \& L. VAINERMAN, Finite Quantum Groupoids and Their
Applications, in New Directions in Hopf Algebras {\it MSRI Publications} {\bf 43} (2002) Cambridge University Press , 211-262.

[NV4] D. NIKSHYCH \& L. VAINERMAN, A Galois correspondence for 
$II_1$-factors and quantum groupoids. {\it JFA} {\bf 178} (2000),
113-142.

[R] RENAULT J.  {\it A groupoid approach to $C^*$-algebras}   Lect.Notes
in Mah. {\bf 793} Springer-Verlag 1980..

[Val0] VALLIN J.M., Unitaire pseudo-multiplicatif associ\'e \`a un 
groupo\" \i de. Applications \`a la moyennabilit\'e. {\it J. of 
Operator Theory} {\bf 44}, No.2 (2000), 347-368.

[Val1] VALLIN J.M. , Groupo\" \i des quantiques finis. {\it Journal 
of Algebra} {\bf 26} (2001), 425-488.

[Val2]  VALLIN J.M. , Multiplicative partial
isometries and finite quantum groupoids. Proceedings of the
Meeting of Theoretical Physicists and Mathematicians,
Strasbourg, 2002. IRMA Lectures in Mathematics and
Theoritical Physics {\bf 2} 189-227. 

[Val3]  VALLIN J.M.  Deformation of finite dimensional quantum  groupoids (math Q.A. 0310265)

[VV] S.VAES \& L.VAINERMAN Extensions of locally compact quantum groups and the bicrossed product construction { \it Advances in Mathematics} {\bf 175} (2003), 1-101.

[Y] YAMANOUCHI T. Duality for actions and co-actions of groupoids on
von-Neumann algebras  {\it Memoirs of the American Mathematical Society}
{\bf 484}(1993).

\end{document}